\documentclass[12pt, a4paper,reqno]{amsart}
\usepackage{amscd,amssymb,amsmath, array,latexsym,amsbsy}

\allowdisplaybreaks[2]
\begin{document}
\def\supp{\operatorname{supp}}
\def\tr{\operatorname{tr}}
\def\lt{\operatorname{{lt}}}
\def\interior{\operatorname{Int}} 
\def\Ind{\operatorname{Ind}} 
\def\Prim{\operatorname{Prim}}
\def\H{\mathcal{H}} 
\def\K{\mathcal{K}} 
\def\N{\mathcal{N}} 
\def\C{\mathbb{C}}
\def\T{\mathbb{T}}
\def\Z{\mathbb{Z}}
\def\R{\mathbb{R}}
\def\P{\mathbb{P}}
\newtheorem{thm}{Theorem}[section]
\newtheorem{cor}[thm]{Corollary}
\newtheorem{prop}[thm]{Proposition}
\newtheorem{lemma}[thm]{Lemma}
\theoremstyle{definition}
\newtheorem{defn}[thm]{Definition}
\newtheorem{remark}[thm]{Remark}
\newtheorem*{remark1}{Remark on hypotheses}
\newtheorem{example}[thm]{Example}
\numberwithin{equation}{section}
\title
[\boldmath Strength of convergence in the orbit space of a transformation
group]
{Strength of convergence in the orbit space of a transformation group}

\author[Archbold]{Robert Archbold}
\address{Department of Mathematical Sciences
\\University of Aberdeen
\\Aberdeen AB24 3UE
\\Scotland
\\United Kingdom
}
\email{r.archbold@maths.abdn.ac.uk}

\author[an Huef]{Astrid an Huef}
\address{School of Mathematics
\\University of New South Wales
\\Sydney, NSW 2052
\\Australia}
\email{astrid@unsw.edu.au} \keywords{Transformation group, orbit
space, $k$-times convergence, spectrum of a $C^*$-algebra,
multiplicity of a representation, trace function}
\subjclass[2000]{22D30, 46L05, 46L30, 46L55, 54H15, 57S05}
\date{August 18, 2005}

\begin{abstract}
Let $(G, X)$ be a  second-countable transformation group with $G$
acting freely on $X$. It is shown that measure-theoretic
accumulation of the action and topological strength of convergence
in the orbit space $X/G$ provide equivalent ways of quantifying
the extent of non-properness of the action. These notions are
linked via the representation theory of the transformation-group
$C^*$-algebra $C_0(X)\times G$.
\end{abstract}

\thanks{This research was supported by grants from the Royal Society
and the University of New South Wales.}

\maketitle

\section{Introduction}
Let $G$ be a locally compact Hausdorff group acting  continuously
on a locally compact Hausdorff space $X$, so that $(G,X)$ is a
transformation group. In \cite{Green2} Green gave a seminal
example of a non-proper action of $G=\R$ on a subspace $X$ of
$\R^3$; this example was generalised in \cite{rieffel} by Rieffel
who replaced the double folding of orbits by  variable folding
using a repetition number associated to each orbit. The main
purpose of this paper is to show that, for free actions, the
fundamental measure-theoretic and topological properties of these
examples persist into the general case. In particular, we show
that the local accumulation of Haar measure (which arises from the
folding of orbits in the case of the Green-Rieffel examples)
corresponds exactly to a topological notion of strength of
convergence in the orbit space $X/G$ from \cite{AD}. The latter
notion is motivated by the phenomenon of strength of convergence
for the Kirillov orbits associated to  a nilpotent Lie group
\cite{ludwig, AKLSS}.

To show that these notions of ``counting by measure'' and
``counting by topology'' give equivalent ways of quantifying the
extent of non-properness of an action, we use the
transformation-group $C^*$-algebra $C_0(X)\times G$. It is
well-known that the representation theory of $C_0(X)\times G$ is
closely related to the properties of the action of $G$ on $X$,
particularly when $G$ and $X$ are second countable.
In \cite[Theorem~3.1]{W} Williams showed that $C_0(X)\times G$ is
liminal (that is,  the image of every irreducible representation
$\pi$ of $C_0(X)\times G$ is the compact operators) if and only if
each orbit is closed in $X$ and all the stability subgroups are
liminal. More generally, Gootman proved in
\cite[Theorem~3.3]{Goot} that $C_0(X)\times G$ is Type \textrm{I}
if and only if each orbit is locally closed in $X$ and each
stability subgroup is Type \textrm{I}. Further, if $C_0(X)\times
G$ is Type \textrm{I} and the action of $G$ on $X$ is free, the
spectrum $(C_0(X)\times G)^\wedge$ of $C_0(X)\times G$ is
homeomorphic to $X/G$ \cite[Lemma~16]{Green2}.

A subset $K$ of $X$ is \emph{wandering} if the set \[\{s\in G:s\cdot K\cap
K\neq\emptyset\}\] is relatively compact in $G$.
The action of $G$ on $X$ is called \emph{proper} if every compact subset of $X$
is wandering, or equivalently,
if the map \[(s,x)\mapsto (s\cdot x, x):G\times X\to X\times X\] is proper in
the sense that inverse images of compact
subsets are compact. If every point of $X$ has a wandering neighbourhood, then
$X$ is a
\emph{Cartan} $G$-space \cite[Definition~1.1.2]{palais}.
For free actions, the conditions on the spectrum of $C_0(X)\times G$
corresponding to proper actions and Cartan $G$-spaces
are that $C_0(X)\times G$ has continuous trace \cite[Theorem~17]{Green2}   and
$C_0(X)\times G$ is a
Fell algebra \cite{aH1},
respectively.
That the action of $G$ on $X$ is proper if and only if the orbit space $X/G$ is
Hausdorff and $X$ is a Cartan $G$-space
\cite[Theorem~1.2.9]{palais} is thus closely related to  a $C^*$-algebra having
continuous trace if and only if it is a
Fell algebra with Hausdorff spectrum \cite[\S4.5]{dixmier}.

More general than a Cartan $G$-space is the notion of an integrable action of
$G$ on $C_0(X)$ from \cite[Definition~1.10]{rieffel}; after translating that
definition to the transformation group $(G, X)$,
an action of $G$ on $X$ is \emph{integrable} if and only if, for each compact
subset $N$ of
$X$,
\begin{equation}
\label{eq-intro}
\sup_{x\in N}\nu\left(\{s\in G:s\cdot x\in N\}\right)<\infty,
\end{equation}
where $\nu$ is a (right) Haar measure on $G$. For free actions,  it
is shown in \cite[Theorem~4.8]{aH2} that the action is integrable if
and only if $C_0(X)\times G$ has bounded trace (that is, there is a
dense ideal $J$ of $C_0(X)\times G$ such that
$\pi\mapsto\tr(\pi(a))$ is bounded on $(C_0(X)\times G)^\wedge$ for
every fixed positive element $a$ of $J$). On the other hand, it
follows from \cite[Theorem~2.6]{ASS} that $C_0(X)\times G$ has
bounded trace if and only if the upper multiplicity $M_U(\pi)$ is
finite for every irreducible representation $\pi$ of $C_0(X)\times
G$ (see \S\ref{sec-prelim} for the definition of the multiplicity
numbers). So, when the action is free and integrable, it is natural
to ask how the finite upper multiplicity numbers of the
representations of $C_0(X)\times G$  are related to the finite
suprema occurring in \eqref{eq-intro}.

Suppose that $(x_n)_{n\geq1}$ is a sequence in $X$
converging to $z\in X$, and that $\Ind\epsilon_{x_n}$  and $\Ind\epsilon_z$
are the corresponding induced representations of $C_0(X)\times G$.
The examples in \cite{Green2, rieffel, deicke, aH2} suggest that the upper and
lower multiplicities $M_U(\Ind\epsilon_z,(\Ind\epsilon_{x_n}))$ and
$M_L(\Ind\epsilon_z,(\Ind\epsilon_{x_n}))$ of $\Ind\epsilon_z$ relative to the
sequence
$(\Ind\epsilon_{x_n})_{n\geq 1}$ depend on the ratios
\begin{equation}\label{eq-intro-ratio}
\frac{\nu(\{s:s\cdot x_n\in
V\})}{\nu(\{s:s\cdot z\in V\})}
\end{equation}
which compare the amount of ``time'' $x_n$ and $z$ spend in a
neighbourhood $V$ of $z$ under the action of $G$.

By contrast, a topological rather than measure-theoretic approach
was adopted in \cite{AD}. Motivated by the notion of strength of
convergence for the Kirillov orbits associated to a nilpotent Lie
group (see \cite[Definition~2.10]{ludwig} and
\cite[Theorem~2.4]{AKLSS}), a sequence $(x_n)_{n\geq 1}$ in $X$ is
said to \emph{converge $k$-times in $X/G$ to $z\in X$} if there
exist $k$ sequences $(t_n^{(1)})_n,(t_n^{(2)})_n,\cdots
,(t_n^{(k)})_n$ in $G$ such that
\begin{enumerate}
\item $t_n^{(i)}\cdot x_n\to z$ as $n\to\infty$ for $1\leq i\leq k$, and
\item if $1\leq i<j\leq k$ then $t_n^{(j)}(t_n^{(i)})^{-1}\to\infty$ as
$n\to\infty$.
\end{enumerate}
It follows from \cite[Theorem~2.3]{AD} that if $(x_n)_{n\geq 1}$
converges $k$-times in $X/G$ to $z$  then
$M_U(\Ind\epsilon_z,(\Ind\epsilon_{x_n}))\geq k$.

The main theorem of this paper shows that, for free actions,  the
measure-theoretic and
topological approaches are  equivalent. Of the five equivalent
conditions of the theorem, only the second involves  the
transformation-group $C^*$-algebra $C_0(X)\times G$; its
representation theory, via the multiplicity numbers, is the key to
passing from topological convergence to measure accumulation.

\begin{thm}\label{thm-main}
Suppose that $(G,X)$ is a free second-countable transformation
group. Let $k$ be a positive integer,  let $z\in X$  and let
$(x_n)_{n\geq 1}$ be a sequence in $X$. Assume that $G\cdot z$ is
locally closed in $X$. Then the following are equivalent:
\begin{enumerate}
\item the sequence $(x_n)_n$ converges $k$-times in $X/G$ to $z$;
\item $M_L(\Ind\epsilon_z,(\Ind\epsilon_{x_n}))\geq k$;
\item for every open  neighbourhood $V$ of $z$  such that
$\{s\in G:s\cdot z\in V\}$  is relatively compact we have
\[
\liminf_n\nu(\{s\in G:s\cdot x_n\in V\})\geq k\nu(\{s\in G:s\cdot z\in V\});
\]
\item there exists a real number $R>k-1$ such that for every open
neighbourhood $V$ of $z$  with $\{s\in G:s\cdot z\in V\}$ relatively
compact we have
\[
\liminf_n\nu(\{s\in G:s\cdot x_n\in V\})\geq R\nu(\{s\in G:s\cdot
z\in V\});
\]
\item there exists a decreasing sequence of basic compact
neighbourhoods $(W_m)_{m\geq 1}$ of $z$ such that, for each
$m\geq1$,
\[
\liminf_n\nu(\{s\in G:s\cdot x_n\in W_m\})> (k-1)\nu(\{s\in G:s\cdot z\in
W_m\}).
\]
\end{enumerate}
\end{thm}
In (3)--(5) above, it is to be understood that the limit infima are
calculated in $[0,\infty]$ and that some of the values $\nu(\{s\in
G:s\cdot x_n\in V\})$ might be infinite. The equivalence of (3),
(4) and (5) demonstrates a surprising emergence of integer values from what
appears in general to be a continuous setting. Note that if the
action is proper, or more generally if $X$ is a Cartan $G$-space,
then the upper multiplicity $M_U(\Ind\epsilon_z)=1$ for all $z\in X$
\cite{Green2, aH1, A}, and hence there can only be $1$-times
convergence in $X/G$ to $z$.

This paper is organised as follows. In \S\ref{sec-prelim} we
establish our conventions and discuss our hypotheses. In
\S\ref{sec-measure} we establish a relationship between the ratios
\eqref{eq-intro-ratio} and the lower and upper multiplicities
$M_L(\Ind\epsilon_z,(\Ind\epsilon_{x_n}))$ and
$M_U(\Ind\epsilon_z,(\Ind\epsilon_{x_n}))$ relative to the
sequence $(\Ind\epsilon_{x_n})_n$ (see Theorems~\ref{thm-M} and
\ref{thm-Mupper}). In \S\ref{sec-top} we link the limit suprema
and infima of measure ratios and topological strength of
convergence. A crucial ingredient  for the results in
\S\ref{sec-measure}  and \S\ref{sec-top} is a new technical lemma
(Lemma~\ref{lem-excision}) which shows that a key property which
we have observed in Green's example \cite{Green2} can be
abstracted to the general situation.

We combine our results in \S\ref{sec-main} to prove
Theorem~\ref{thm-main}. The first corollary of
Theorem~\ref{thm-main} is the analogous result involving $k$-times
convergence of a subsequence of $(x_n)_n$, the upper multiplicity
$M_U(\Ind\epsilon_z,(\Ind\epsilon_{x_n}))$  and limit suprema of
the measure ratios (Theorem~\ref{thm-main-upper}).

We apply our results to examples in \S\ref{sec-examples}. In
particular, we consider Rieffel's example of a free transformation
group  \cite[Example 1.18]{rieffel} where the orbit space $X/G$
consists of a sequence  $(G\cdot x_n)_{n\geq 1}$ converging to
$G\cdot x_0$.  For each $n\geq 1$ the orbit has $L_n+1$ folds, where
$L_n$ is a repetition number. We show that
$M_L(\Ind\epsilon_{x_0},(\Ind\epsilon_{x_n}))=\liminf{(L_n+1)}$ and
$M_U(\Ind\epsilon_{x_0},(\Ind\epsilon_{x_n}))=\limsup{(L_n+1)}$.

In Appendix~\ref{appendix}
we establish sequence versions of
\cite[Propositions 2.2 and 2.3]{AS} which are needed in \S\ref{sec-measure} and
in \S\ref{sec-main}.

\section{Preliminaries}\label{sec-prelim}
Throughout,  $(G,X)$ is a locally compact Hausdorff transformation
group: thus $G$ is a locally compact Hausdorff group and $X$ is a
locally compact Hausdorff space together with a jointly continuous
map $(s,x)\mapsto s\cdot x$ from $G\times X$ to $X$ such that
$s\cdot (t\cdot x)=st\cdot x$ and $e\cdot x=x$.  In all our main results the
action is assumed to be \emph{free}, that is  $s\cdot x=x$ implies $s=e$. We
also assume that both $G$
and $X$ are second countable, and hence both are normal spaces
\cite[Chapter 4, Lemma~1]{K}.

The action of $G$ on $X$ lifts to a strongly continuous action
$\lt$ of $G$ by automorphisms of $C_0(X)$ given by
$\lt_s(f)(x)=f(s^{-1}\cdot x)$.
The associated transformation-group $C^*$-algebra $C_0(X)\times G$ is
the $C^*$-algebra which is universal for the covariant representations
of the $C^*$-dynamical system $(C_0(X), G, \lt)$.
More concretely, $C_0(X)\times G$ is the
enveloping $C^*$-algebra of the Banach $*$-algebra $L^1(G, C_0(X))$ of
functions $f:G\to C_0(X)$ which are integrable with respect to a
fixed left Haar measure $\mu$ on $G$, with
multiplication and involution  given by
\[
f*g(s):=\int_G f(r)\lt_r(g(r^{-1}s))\ d\mu(r)\quad\text{and}
\quad f^*(s):=\Delta(s^{-1})\lt_s(f(s^{-1}))^*,
\]
where $\Delta$ is the modular function associated with $\mu$.

Suppose that $G$ acts freely on $X$ and let $\nu$ be the right Haar measure such that $\nu(E)=\mu(E^{-1})$.
As in \cite{W, W2, aH2}, for each $x\in X$ we realise the induced
representation $\Ind\epsilon_x$ (where $\epsilon_x:C_0(X)\to\C$ is
evaluation at $x$) on the Hilbert space $L^2(G,\nu)$ as the
integrated form
\[\Ind\epsilon_x=\tilde\epsilon_x\times \lambda,\]
where $(\tilde\epsilon_x(f)\xi)(s)=f(s\cdot x)\xi(s)$ and
$(\lambda_t\xi)(s)=\Delta(t)^{1/2}\xi(t^{-1}s)$ for $f\in C_0(X)$,
$\xi\in L^2(G,\nu)$ and $r,s\in G$ (see \cite[Lemma~4.14]{W}). Thus
\begin{align}\label{eq-pix}
(\Ind\epsilon_x(f)\xi)(s)&=\int_G f(t,s\cdot
x)\xi(t^{-1}s)\Delta(t)^{1/2}\,
d\mu(t)\notag\\
&=\int_G f(sw^{-1}, s\cdot x)\Delta(sw^{-1})^{1/2}\xi(w)\, d\nu(w)
\end{align}
for $f\in C_c(G,X)$ and $\xi\in L^2(G,\nu)$.  Note that
$\Ind\epsilon_x$ is irreducible by, for example,
\cite[Proposition~4.2]{W}.

 Let $A$ be a $C^*$-algebra and let $\pi$ be
an irreducible representation of $A$; we use
the same symbol $\pi$ to denote the
unitary equivalence class of $\pi$ in the spectrum $\hat A$ of $A$.
If $\pi_1$ and $\pi_2$ are equivalent
irreducible representations then
$\tr(\pi_1(a))=\tr(\pi_2(a))$ for all $a$ in the postitive cone $A^+$ of $A$,
and so we may
write unambiguously the expression $\tr(\pi(a))$ whenever $\pi\in
{\hat A}$ and $a\in A^+$.

We now recall the definitions of upper and lower multiplicity from
\cite{A}, \cite{AS}. Unless stated otherwise, we shall always regard the
Banach dual $A^*$ as being equipped with the weak*-topology. We
denote by ${\mathcal N}$ the weak$^*$-neighbourhood base at zero
in $A^*$ consisting of all open sets of the form
$$N=\{\psi\in A^*\>:\>\vert\psi(a_i)\vert<\epsilon,\>1\leq i\leq n\}$$
where $\epsilon>0$ and $a_1,\ldots,a_n \in A$. Let $P(A)$ be the set
of pure states of $A$ and let $\theta:\>P(A)\to {\hat A}$ be the
continuous, open mapping given by $\theta(\phi)=\pi_{\phi}$ where
$\pi_{\phi}$ is the  GNS representation associated with $\phi$
\cite[3.4.11]{dixmier}.

We write $\P=\mathbb{N}\setminus\{0\}$. Let $\pi\in {\hat A}$. First of all, we
give descriptions of the
upper and lower multiplicities $M_U(\pi)$ and $M_L(\pi)$. Let
$\phi$ be a pure state of $A$ associated with $\pi$ and let
$N\in{\mathcal N}$. Let
$$V(\phi,N) = \theta((\phi+N)\cap P(A)),$$ an open neighbourhood
of $\pi$ in ${\hat A}$. For $\sigma\in {\hat A}$ let $\H_\sigma$ be the Hilbert
space of $\sigma$, and let
$${\rm Vec}(\sigma,\phi,N) = \{\eta\in H_{\sigma}\>:\>\Vert\eta\Vert=1,\>
\langle \sigma(\cdot)\eta,\eta\rangle \in \phi+N\}.$$ Note that
${\rm Vec}(\sigma,\phi,N)$ is nonempty if and only if $\sigma\in
V(\phi,N)$. For $\sigma\in V(\phi,N)$ we define $d(\sigma,\phi,N)$
to be the supremum (in $\P\cup\{\infty \}$) of the cardinalities
of finite orthonormal subsets of ${\rm Vec}(\sigma,\phi,N)$. It is
convenient to define $d(\sigma,\phi,N)=0$ for $\sigma\in {\hat
A}\setminus V(\phi,N)$.

>>From \cite[\S2 and Proposition~3.4]{A} we have
$$M_U(\pi)=\inf_{N\in\N}\left(\limsup_{\sigma\to\pi}\>d(\sigma,\phi,N)\right)\in
\P\cup\{\infty \}$$ and, if $\pi$ is not open in ${\hat A}$,
$$M_L(\pi)=\inf_{N\in\N}\left(\liminf_{\sigma\to\pi,\>\sigma\neq\pi}
d(\sigma,\phi,N)\right)\in \P\cup\{\infty\}.$$ As noted in
\cite[Lemma~2.1]{A}, $M_U(\pi)$ and $M_L(\pi)$ are independent of
the choice of $\phi$. Elementary examples which motivate these definitions and
illustrate the computations are given in \cite[$\S2$]{A}.

Now suppose, in addition, that
$\Omega=(\pi_{\alpha})_{\alpha\in\Lambda}$ is a net in ${\hat A}$.
For $N\in\N$ let
$$M_U(\phi,N,\Omega)=\limsup_{\alpha}\>d(\pi_{\alpha},\phi,N)\in
{\mathbb N}\cup\{\infty\}.$$ Note that if $N'\in {\N}$ and
$N'\subset N$ then $M_U(\phi,N',\Omega)\leq M_U(\phi,N,\Omega)$.
We define
$$M_U(\pi,\Omega)=\inf_{N\in\N}M_U(\phi,N,\Omega)\in {\mathbb N}\cup\{\infty\}$$
(which is independent of the choice of $\phi$ by an argument similar
to that used in the proof of \cite[Lemma~2.1]{A}). Similarly, for
$N\in\N$, let
$$M_L(\phi,N,\Omega)=\liminf_{\alpha}\>d(\pi_{\alpha},\phi,N)\in
{\mathbb N}\cup\{\infty\}.$$ Then $M_L(\phi,N,\Omega)$ decreases
with $N$ and we define
$$M_L(\pi,\Omega)=\inf_{N\in\N}M_L(\phi,N,\Omega)\in{\mathbb N}\cup\{\infty\}$$
(which is, again, independent of the choice of $\phi$). Note that
it is not required that $\Omega$ converge to $\pi$.  However it
follows from these definitions that $M_U(\pi,\Omega)>0$ if and
only if $\pi$ is a cluster point of $\Omega$, and that
$M_L(\pi,\Omega)>0$ if and only if $\Omega$ converges to $\pi$.

We have the inequalities
$$M_L(\pi,\Omega)\leq M_U(\pi,\Omega)\leq M_U(\pi)$$
and, if $\Omega$ is convergent to $\pi$ but eventually
$\pi_\alpha\not=\pi$,
$$M_L(\pi)\leq M_L(\pi,\Omega)$$
(see \cite[Proposition~2.1]{AS}). Also, if $\Omega_0$ is a subnet of
$\Omega$ then
$$M_L(\pi,\Omega)\leq M_L(\pi,\Omega_0)\leq M_U(\pi,\Omega_0)\leq
M_U(\pi,\Omega).$$

\subsection*{Remark on hypotheses}
In our major results, we assume that  the
action of $G$ on $X$ is free. Typically we focus on a fixed $z\in X$
together with a sequence $(x_n)_{n\geq 1}$ in $X$ such that  the
orbit $G\cdot z:=\{s\cdot z:s\in G\}$ is locally closed in $X$  (in
the sense that $G\cdot z$ is relatively open in its closure) and
$G\cdot x_n\to G\cdot z$ in the orbit space $X/G$.

Lemma~\ref{lem-remark} explains the hypotheses in our main
Theorem~\ref{thm-main}. In particular, if $G\cdot z$ is not
locally closed then one cannot form the measure ratios (1.2). For
each $x\in X$, let $\phi_x:G\to G\cdot x$ be the map $s\mapsto
s\cdot x$.

\begin{lemma}\label{lem-remark}
Let $(G,X)$ be a second-countable transformation group and $\nu$ a
right Haar measure on $G$. Let $z\in X$ and suppose that the
stability subgroup $S_z:=\{s\in G:s\cdot z=z\}$ is compact.  Then
the following are equivalent:
\begin{enumerate}
\item the orbit $G\cdot z$ is not locally closed in $X$;
\item for every $k\in\P$, the sequence $z,z,z,\dots$ converges
$k$-times in $X/G$ to $z$;
\item for every open neighbourhood $V$ of $z$,
$\nu(\phi_z^{-1}(V))=\infty$;
\item for every open neighbourhood $V$ of $z$,
$\phi_z^{-1}(V)$ is not relatively compact in $G$.
\end{enumerate}
\end{lemma}
\begin{proof}
Let $(V_n)_{n\geq 1}$ be a decreasing sequence of basic open
neighbourhoods of $z$ in $X$ and let $(K_n)_{n\geq 1}$ be an
increasing sequence of compact subsets of $G$ such that
$G=\cup_{n\geq 1}\interior(K_n)$.

(1) $\Longrightarrow$ (2). Suppose that $G\cdot z$ is not locally
closed. Then $W\cap (\overline{G\cdot z}\setminus G\cdot
z)\neq\emptyset$ for every neighbourhood $W$ of $z$. Let $k\geq1$.
We will construct $k$ sequences $(t_n^{(i)})_{n\geq 1}\subset G$,
where $1\leq i\leq k$, such that $t_n^{(i)}\cdot z\in V_n$ and
$t_n^{(j)}(t_n^{(i)})^{-1}\notin K_n$ for each $n\geq 1$.

Temporarily fixing $n$, we construct $t_n^{(i)}$ as follows.  Let
$t_n^{(1)}=e$. Since $G\cdot z$ is not locally closed there exists
$y\in V_n\cap (\overline{G\cdot z}\setminus G\cdot z)$. Since $y$
is in the closure of $G\cdot z$ and $V_n$ is open, given any
compact subset $K$ of $G$ there exists $t_K\in G\setminus K$ such
that $t_K\cdot z\in V_n$.  So there exists $t_n^{(2)}\in
G\setminus K_n$ such that $t_n^{(2)}\cdot z\in V_n$. Proceeding
inductively we obtain $t_n^{(2)},t_n^{(3)} \dots, t_n^{(k)}$ such
that
\[t_n^{(j)}\cdot z\in V_n\quad\text{and}\quad t_n^{(j)}\in G\setminus \left(
\cup_{i=1}^{j-1} K_nt_n^{(i)} \right)
\]
for $2\leq j\leq k$.  Letting $n$ run, it follows easily from the
properties of $(V_n)_{n\geq 1}$ and $(K_n)_{n\geq 1}$ that
$z,z,z,\dots$ converges $k$-times in $X/G$ to $z$.

(2) $\Longrightarrow$ (3). Suppose that (2) holds. Let $V$ be an
open neighbourhood of $z$ and $M>0$.  By the continuity of the
action on the locally compact Hausdorff space $X$, there exists an
open neighbourhood $U$ of $z$ and a compact neighbourhood $K$ of $e$
in $G$ such that $K\cdot U\subset V$. Choose $k\in\P$ such that
$k\nu(K)>M$. By (2) there exist $k$ sequences $(t_n^{(i)})_{n\geq 1}$
 such that $t_n^{(i)}\cdot z\to z$ as $n\to\infty$ for each $1\leq i\leq k$, and
\[ t_n^{(j)}(t_n^{(i)})^{-1}\to\infty\text{\ as\ }
n\to\infty\quad (1\leq i<j\leq k).
\]
Hence there exists $n_0$ such that  $t_{n_0}^{(i)}\cdot z\in U$ for
$1\leq i\leq k$ and $t_{n_0}^{(j)}(t_{n_0}^{(i)})^{-1}\in
G\setminus(K^{-1}K)$ for $1\leq i<j\leq k$. Then
$Kt_{n_0}^{(i)}\cdot z\subset K\cdot U\subset V$ for $1\leq i\leq k$
and $Kt_{n_0}^{(i)}\cap  Kt_{n_0}^{(j)}=\emptyset$ unless $i=j$, and
hence $\nu(\phi_z^{-1}(V))\geq k\nu(K)>M$.  Since $M$ was arbitrary,
the result follows.

(3) $\Longrightarrow$ (4). Compact subsets have finite Haar measure,
so this is immediate.

 (4) $\Longrightarrow$ (1). Suppose that $G\cdot z$ is locally closed
 in $X$. Then $G\cdot z$ is a relatively open subset of the locally
 compact space $\overline{G\cdot z}$ and hence $G\cdot z$ is locally
 compact.  Thus $(G,G\cdot z)$ is a second-countable
 locally compact Hausdorff transformation group.  In particular, it
 follows from \cite[Theorem~1]{Gli} that the map $\psi_z:G/S_z\to
 G\cdot z:sS_z\mapsto s\cdot z$ is a homeomorphism.
Let $U$ be an open subset  of $X$ such that $U\cap\overline{G\cdot
z}=G\cdot z$. Let $N$ be a compact neighbourhood of $z$ in $X$ such
that $N\subset U$. Then $N\cap \overline{G\cdot z}=N\cap G\cdot z$
is a compact subset of $G\cdot z$. Hence $\psi_z^{-1}(N)$ is compact
in $G/S_z$. Since $S_z$ is compact, the quotient map $q_z:G\to G/S_z$ is
proper, and hence $\phi_z^{-1}(V)=q_z^{-1}(\psi_z^{-1}(V))$ is compact in $X$. Let
$V=\interior N$; then $\phi_z^{-1}(V)$ is relatively compact in $G$.
\end{proof}

If $M_U(\Ind\epsilon_z)$ were finite, then we could assume,
without loss of generality, that all orbits are closed in $X$. To
see this, note that $\Ind\epsilon_z$ restricts to an irreducible
representation of the bounded-trace ideal $J$ of $C_0(X)\times G$
\cite[Theorem~2.8]{ASS}, and $J$ is canonically isomorphic to
$C_0(Y)\times G$ for some $G$-invariant open subset $Y$ of $X$
\cite[Theorem~5.8]{aH2}. Since $J$ (and hence $C_0(Y)\times G$) is
liminal, the orbits in $Y$ are relatively closed in $Y$. By
replacing $X$ with $Y$ we are then in the  situation of closed
orbits without changing any dynamical or representational
properties at or near $z$; in particular, the multiplicities for
$\Ind\epsilon_z$ are the same whether we compute them in the ideal
$C_0(Y)\times G$ or in $C_0(X)\times G$
\cite[Proposition~5.3]{ASS}. However, we are also interested in
examples where $M_U(\Ind\epsilon_z)=\infty$ but the lower
multiplicity $M_L(\Ind\epsilon_z)$ is finite, and in examples
where $M_L(\Ind\epsilon_z)=M_U(\Ind\epsilon_z)=\infty$.



\section{Measure ratios and bounds on multiplicity numbers}
\label{sec-measure}

We will frequently use the following well-known properties of a
locally compact Hausdorff space X (see, for example,
\cite[Chapter~5, Theorem~18]{K}): if U is a neighbourhood of a
compact subset $N$ of $X$, then $U$ contains a compact neighbourhood
$V$ of $N$ and there is a continuous function $f:X\to [0,1]$ which
is $1$ on $N$ and $0$ on $X\setminus V$ (and hence on $X\setminus
U$).

For $x\in X$ and $f\in C_0(X)$ we
define the function $f_x:G\to \C$ by
 $f_x(s)=(f\circ\phi_x)(s)=f(s\cdot x)$.

\begin{thm}\label{thm-Msquared}
Suppose that $(G,X)$ is a free second-countable transformation
group. Let  $M\in\R$ with $M\geq 1$, let $z\in X$  and let
$(x_n)_{n\geq 1}$ be a sequence in $X$. Assume that $G\cdot z$ is
locally closed in $X$.  Suppose that there exists an open
neighbourhood $V$ of $z$ in $X$ such that $\phi_z^{-1}(V)$ is
relatively compact and
\begin{equation*}
\nu(\phi_{x_n}^{-1}(V))\leq M\nu(\phi_z^{-1}(V))
\end{equation*}
frequently. Then $M_L(\Ind\epsilon_z,(\Ind\epsilon_{x_n}))\leq\lfloor M^2\rfloor $.
\end{thm}

\begin{proof}
Fix $\epsilon>0$ such that $M^2(1+\epsilon)^2<\lfloor M^2\rfloor +1$.
We will build a function $D\in C_c(G\times X)$
such that $\Ind\epsilon_z(D^**D)$ is a rank-one projection and
\[
\tr(\Ind\epsilon_{x_n}(D^**D))<M^2(1+\epsilon)^2<\lfloor M^2\rfloor +1
\] frequently.
(The function $D$ is similar to the ones used in \cite[Proposition~4.2]{W2} and \cite[Proposition~4.5]{aH2}.)
By the generalised lower semi-continuity result of \cite[Theorem~4.3]{AS} we will have
\begin{align*}
\liminf_n \tr(\Ind\epsilon_{x_n}(D^**D))&\geq
M_L(\Ind\epsilon_z,(\Ind\epsilon_{x_n}))\tr(\Ind\epsilon_z(D^**D))\\&=
M_L(\Ind\epsilon_z,(\Ind\epsilon_{x_n})),
\end{align*}
and the theorem will follow.

Let $\delta>0$ such that
\begin{equation*}
\delta<\frac{\epsilon\nu(\phi_z^{-1}(V))}{1+\epsilon}<\nu(\phi_z^{-1}(V)).
\end{equation*}
By the regularity of the measure $\nu$ there exists a compact subset $W$ of the
open set $\phi_z^{-1}(V)$ such that
\begin{equation*}
0<\nu(\phi_z^{-1}(V))-\delta<\nu(W).
\end{equation*}
Since $W$ is compact, there is a compact neighbourhood $W_1$ of
$W$ contained in $\phi_z^{-1}(V)$ and a continuous function
$g:G\to[0,1]$ such that $g$ is identically one on $W$ and is
identically zero off the interior of $W_1$. Then
\begin{equation*}
\nu(\phi_z^{-1}(V))-\delta<\nu(W)\leq \int_G g(t)^2\, dt=\|g\|_2^2,
\end{equation*}
and hence
\begin{equation}\label{eq-estimate}
\frac{\nu(\phi_z^{-1}(V))}{\|g\|_2^2}<
1+\frac{\delta}{\|g\|_2^2}<1+\frac{\delta}{\nu(\phi_z^{-1}(V))-\delta}<1+\epsilon.
\end{equation}
Since $G\cdot z$ is locally closed in $X$ it follows from
\cite[Theorem~1]{Gli}, applied to the locally compact Hausdorff
transformation group $(G, G\cdot z)$, that $\phi_z$ is a
homeomorphism of $G$ onto $G\cdot z$. So there is a continuous
function $g_1:W_1\cdot z\to[0,1]$ such that $g_1(t\cdot z)=g(t)$ for
$t\in W_1$. Since $W_1\cdot z$ is a compact subset of the locally
compact Hausdorff space $X$, it follows from Tietze's Extension
Theorem (applied to the one-point compactification of $X$ if
necessary) that $g_1$  can be extended to a continuous function
$g_2:X\to [0,1]$. Because $W_1\cdot z$ is a compact subset of the
open set $V$, there exists a compact neighbourhood $P$ of $W_1\cdot
z$ contained in $V$ and a continuous function $h:X\to [0,1]$ such
that $h$ is identically one on $W_1\cdot z$ and is identically zero
off the interior of $P$.
 Note that $h$ has compact support contained in $P$. We set
\[
f(x)=h(x)g_2(x).
\]
Then $f\in C_c(X)$ with $0\leq f\leq 1$ and  $\supp f\subset\supp
h\subset P\subset V$. Note that
\begin{equation}\label{eq-reciprocal}
 \|f_z\|^2_2=\int_G f(t\cdot z)^2\, dt=\int_G h(t\cdot z)^2g_2(t\cdot z)^2\, dt\geq\int_{W_1} g(t)^2\, dt =\|g\|^2_2
\end{equation}
since $h$ is identically one on $W_1\cdot z$ and $g$ has support
inside $W_1$. We now set
\[
F(x)=\frac{f(x)}{\|f_z\|_2}.
\]
Now $F\in C_c(X)$, $\|F_z\|_2=1$ and $F_x(s)=F(s\cdot x)\neq 0$
implies $s\in\phi_x^{-1}(V)$ by our choice of $h$. Since
$\phi_z^{-1}(V)$ is relatively compact, $\supp F_z$ is compact.

Choose $b\in C_c(G\times X)$ such that $0\leq b\leq 1$ and $b$ is identically
one on the set $(\supp F_z)(\supp F_z)^{-1}\times\supp F$.
Set
\[
B(t,x)=F(x)F(t^{-1}\cdot x)b(t^{-1},x)\Delta(t)^{-1/2}\quad \text{and}\quad
D=\frac{1}{2}(B+B^*).
\]
Using \eqref{eq-pix}, we have
\[
(\Ind\epsilon_x(B)\xi)(s) =\int_G F(s\cdot x)F(w\cdot x)b(ws^{-1},s\cdot
x)\xi(w)\, d\nu(w)\] for $\xi\in L^2(G,\nu)$. So
\[
(\Ind\epsilon_x(D)\xi)(s)=\frac{1}{2}F(s\cdot x)\int_G F(w\cdot x)\big(b(ws^{-1},s\cdot
x)+b(sw^{-1}, w\cdot x)\big) \xi(w)\, d\nu(w).
\]

If $s,w\in \supp F_z$ then $s\cdot z,w\cdot z\in\supp F$ by the continuity of
the action, and hence $b(ws^{-1},s\cdot z)+b(sw^{-1}, w\cdot z)=2$.  It follows
that $\Ind\epsilon_z(D)(\xi)=(\xi, F_z)F_z$.  Thus $\Ind\epsilon_z(D)$, and
hence $\Ind\epsilon_z(D^**D)$, is the rank-one projection determined by the unit vector
$F_z\in L^2(G,\nu)$.

Choose a subsequence $(x_{n_i})_i$ of $(x_n)_n$ such that
\begin{equation*}
\nu(\phi_{x_{n_i}}^{-1}(V))\leq M\nu(\phi_z^{-1}(V))
\end{equation*}
for all $i\geq 1$ and set
$
E_i=\{s\in G: F(s\cdot x_{n_i})\neq 0\}
$.
Then each $E_i$ is open (hence measurable) with
\[\nu(E_i)\leq\nu(\phi_{x_{n_i}}^{-1}(V))\leq M\nu(\phi_z^{-1}(V))<\infty\]
and
\begin{equation}\label{eq-moved}
\int_G F(s\cdot x_{n_i})^2\, d\nu(s)\leq\frac{\nu(E_i)}{\|f_z\|^2_2}
\leq\frac{M\nu(\phi_z^{-1}(V))}{\|g\|^2_2}
\end{equation}
using \eqref{eq-reciprocal}. The function
\[
S_i(s, w):=\frac{1}{2}F(s\cdot x_{n_i}) F(w\cdot x_{n_i})\big(b(ws^{-1},s\cdot x_{n_i})+b(sw^{-1},
w\cdot x_{n_i})\big)
\]
is continuous  and is bounded by $\|F\|_\infty^2$ because $0\leq b\leq 1$.
Thus  $S_i\in L^2(G\times G)$ and  $\Ind\epsilon_{x_{n_i}}(D)$ is the self-adjoint Hilbert-Schmidt operator on $L^2(G,\nu)$ with kernel $S_i$.
It follows that $\Ind\epsilon_{x_{n_i}}(D^**D)$ is a trace-class operator with
\[\tr(\Ind\epsilon_{x_{n_i}}(D^**D))=\|S_i\|_2^2\]
(see, for example, \cite[Proposition~3.4.16]{ped}).
An application of Fubini's Theorem gives
\begin{align}
\tr(&\Ind\epsilon_{x_{n_i}}(D^**D))\label{eq-traceD}\\
&= \frac{1}{4}\int\int F(s\cdot x_{n_i})^2 F(w\cdot x_{n_i})^2\big(
b(ws^{-1},s\cdot x_{n_i})+b(sw^{-1},w\cdot x_{n_i})  \big)^2\, d\nu(w)\,
d\nu(s)\notag\\
&\leq \Big(\int_G F(s\cdot x_{n_i})^2\, d\nu(s)\Big)^2\notag\\
&\leq\frac{M^2\nu(\phi_z^{-1}(V))^2}{\|g\|_2^4}\quad\quad\quad\text{(using \eqref{eq-moved})}\notag\\
&<M^2(1+\epsilon)^2 \quad\quad\quad\quad\text{(using \eqref{eq-estimate}).}\notag
\end{align}
Now
\begin{equation*}
M_L(\Ind\epsilon_z,(\Ind\epsilon_{x_n}))\leq \liminf_n \tr(\Ind\epsilon_{x_n}(D^**D))\leq
M^2(1+\epsilon)^2<\lfloor M^2\rfloor+1,
\end{equation*}
and hence $M_L(\Ind\epsilon_z,(\Ind\epsilon_{x_n}))\leq\lfloor M^2\rfloor$.
\end{proof}


Our aim is to use Theorem~\ref{thm-Msquared} to bootstrap to a
result where  the relative lower multiplicity is bounded by $\lfloor
M\rfloor$ rather than $\lfloor M^2\rfloor$ (see Theorem~\ref{thm-M}).
To do this we will need
to be able to cut away any limits of $(\Ind\epsilon_{x_n})_n$ that are
different from $\Ind\epsilon_z$ so that we can
apply the following crucial Lemma~\ref{lem-excision}.


\begin{lemma}\label{lem-excision}
Suppose that $(G,X)$ is a  transformation group. Let $W$ be a compact neighbourhood of $z\in X$, $K$ a compact subset
of $G$, and $U$ an open neighbourhood of $\phi_z^{-1}(W)$ in $G$.
Let $(x_n)_{n\geq 1}$ be a sequence in $X$ such  that $G\cdot x_n\to
G\cdot z$ and $G\cdot z$ is the unique limit of $(G\cdot x_n)_{n\geq 1}$ in
$X/G$. There exists $n_0$ such that  for every $n\geq n_0$ and every $s\in\phi_{x_n}^{-1}(W)$ there exists $r\in\phi_z^{-1}(W)$ such that
$Ks\cap\phi_{x_n}^{-1}(W)\subset Ur^{-1}s$.
\end{lemma}

\begin{proof}
Suppose that there is no such $n_0$.  Then, by passing to a subsequence, for
each $n$ there exists $s_n\in\phi_{x_n}^{-1}(W)$ such that
\begin{equation}\label{eq-excision}
Ks_n\cap\phi_{x_n}^{-1}(W)\not\subset Ur^{-1}s_n
\end{equation}
for all $r\in\phi_z^{-1}(W)$.

Note that $s_n\cdot x_n\in W$, so by passing to a further subsequence we may
assume that $s_n\cdot x_n$ converges to some $y\in W$.  Since $G\cdot z$ is the
unique limit point of $(G\cdot x_n)$ in  $X/G$, we have $y=r\cdot z$ for some
$r\in G$.  Since $y\in W$ we have $r\in\phi_z^{-1}(W)$, and for this $r$
Equation~\ref{eq-excision} holds for all $n$. So there exist $k_n\in K$ such
that $k_ns_n\in\phi_{x_n}^{-1}(W)$ but $k_n\notin Ur^{-1}$.  By passing to a
further subsequence $(k_n)_n$ converges to $k\in K$.  Then $k_n s_n\cdot x_n\to
k\cdot y=kr\cdot z\in W$. But now $k\in\phi_z^{-1}(W)r^{-1}\subset Ur^{-1}$ and
$k_n\notin Ur^{-1}$ which is impossible since $U$ is open and $k_n\to k$.
\end{proof}


In the next Lemma we consider a potentially bad neighbourhood
$V$ of $z$, where the measure of $\phi_z^{-1}(\overline{ V})$
might be much larger than that of $\phi_z^{-1}(V)$, and show that
we can at least find a nicer, controlled, neighbourhood contained
in it.


\begin{lemma}\label{lem-badV}
Suppose that $(G,X)$ is a  transformation group. Let $\gamma>0$,
$z\in Z$ and let $V$ be an open neighbourhood of $z$ in $X$ such
that $\nu(\phi_z^{-1}(V))<\infty$. Then there exists an open
relatively compact neighbourhood $V_1$ of $z$ such that
$\overline{V_1}\subset V$ and
\begin{equation}\label{eq-lem-gamma}
\nu(\phi_z^{-1}(V))-\gamma
<\nu(\phi_z^{-1}(V_1))
\leq\nu(\phi_z^{-1}(\overline{V_1}))
\leq\nu(\phi_z^{-1}(V))
<\nu(\phi_z^{-1}(V_1))+\gamma.
\end{equation}
\end{lemma}

\begin{proof}
By the regularity of the measure $\nu$ we can choose a compact
subset $W$ of $\phi_z^{-1}(V)$ such that $e\in W$ and
$\nu(W)>\nu(\phi_z^{-1}(V))-\gamma$.  Note that $W\cdot z$ is a
compact subset of $G\cdot z$ and hence of $X$. Hence there exists an
open relatively compact neighbourhood $V_1$ of $W\cdot z$ such that
$\overline{V_1}\subset V$. Now
\[\nu(\phi_z^{-1}(V))-\gamma<\nu(W)\leq\nu(\phi_z^{-1}(V_1))
\leq\nu(\phi_z^{-1}(\overline{V_1})) \leq\nu(\phi_z^{-1}(V))\] and
\eqref{eq-lem-gamma} follows.
\end{proof}


The following result is a sharpening of \cite[Proposition 2.1(i)]{ASS}.
We will use it to cut away multiple limits of our sequences.


\begin{prop}\label{prop-cutaway}
Let $A$ be a $C^*$-algebra and let $\pi\in\hat A$.  Suppose that
$\Omega=(\pi_{\alpha})_{\alpha\in\Lambda}$ is a net in $\hat A$
which is convergent to $\pi$ and that $M_L(\pi,\Omega)=k<\infty$.
Then $\{\pi\}$ is open in the set $L(\Omega)$ of limits of
$\Omega$.
\end{prop}

\begin{proof}
The direct way to see this is to re-work the proof of
\cite[Proposition 2.1(i)]{ASS}.
An alternative shorter argument is as follows. By \cite[Proposition
2.3]{AS}, there is a subnet $\Omega_1$ of $\Omega$ such that
$$M_U(\pi,\Omega_1)=M_L(\pi,\Omega)<\infty.$$ By \cite[Proposition
2.1(i)]{ASS}, $\{\pi\}$ is open in $L(\Omega_1)$. But
$L(\Omega)\subset L(\Omega_1)$ and so $\{\pi\}$ is open in
$L(\Omega)$.
\end{proof}

Theorem~\ref{thm-M} has the same hypotheses as
Theorem~\ref{thm-Msquared} but a stronger conclusion;  the
proof requires Theorem~\ref{thm-Msquared} and some estimates based on Lemmas~\ref{lem-excision} and \ref{lem-badV}.


\begin{thm}\label{thm-M}
Suppose that $(G,X)$ is a free second-countable transformation
group. Let $M\in\R$ with $M\geq 1$, let $z\in X$ and let
$(x_n)_{n\geq 1}$ be a sequence in $X$. Assume that $G\cdot z$ is
locally closed in $X$. Suppose that there exists an open
neighbourhood $V$ of $z$ in $X$ such that $\phi_z^{-1}(V)$ is
relatively compact and
\begin{equation*}
\nu(\phi_{x_n}^{-1}(V))\leq M\nu(\phi_z^{-1}(V))
\end{equation*}
frequently. Then $M_L(\Ind\epsilon_z,(\Ind\epsilon_{x_n}))\leq\lfloor M\rfloor$.
\end{thm}

\begin{proof}
If $\Ind\epsilon_{x_n}\not\to\Ind\epsilon_z$ then $M_L(\Ind\epsilon_z,(\Ind\epsilon_{x_n}))=0<\lfloor M\rfloor$. So we assume from now on that $\Ind\epsilon_{x_n}\to\Ind\epsilon_z$.

Our first claim is that $G\cdot x_n\to G\cdot z$. To see this, suppose that $G\cdot x_n\not\to G\cdot z$.  Then there exists an open neighbourhood $U_0$ of $G\cdot z$ such that $G\cdot x_n$ is frequently not in $U_0$.  Let $q:X\to X/G$ be the quotient map. Then $U_1=q^{-1}(U_0)$ is an open $G$-invariant neighbourhood of $z$ and  $x_n\notin U_1$ frequently.  Note that $C_0(U_1)\times G$ is isomorphic to a closed two-sided ideal $I$ of $C_0(X)\times G$ and  $I\subset\ker(\Ind\epsilon_{x_n})$ whenever $x_n\notin U_1$. Hence $\Ind\epsilon_{x_n}\notin\hat I$ frequently.  But $\hat I$ is an open neighbourhood of $\Ind\epsilon_z$, so $\Ind\epsilon_{x_n}\not\to\Ind\epsilon_z$.

Next we claim  that we may assume, without loss of generality, that $G\cdot z$ is the unique limit of $(G\cdot x_n)_n$ in $X/G$.  To see this, note that $M_L(\Ind\epsilon_z,(\Ind\epsilon_{x_n}))\leq \lfloor M^2\rfloor<\infty$
by Theorem~\ref{thm-Msquared}.
Hence, by Proposition~\ref{prop-cutaway}, $\{\Ind\epsilon_z\}$ is open in
the set of limits of $(\Ind\epsilon_{x_n})_n$. So there is an open
neighbourhood $U_2$ of $\Ind\epsilon_z$ in the spectrum of $C_0(X)\times G$
such that $\Ind\epsilon_z$ is the unique limit of $(\Ind\epsilon_{x_n})_n$ in $U_2$.

Write $(X/G)^\sim$ for the $\mathrm{T}_0$-isation of $X/G$, so that $(X/G)^\sim$ is the quotient of $X$ obtained by identifying points of $X$ with equal orbit closures.  It follows, for example from \cite[Lemmas~4.5 and 4.10]{W} that the map
\[
(X/G)^\sim\to\Prim(C_0(X)\times G)
\]
given by $[x]\mapsto\ker(\Ind\epsilon_x)$
is a homeomorphism.  Since the topology on $(C_0(X)\times G)^\wedge$ is the topology pulled back from $\Prim(C_0(X)\times G)$, it is straightforward to see that $\Ind:X\to (C_0(X)\times G)^\wedge:x\mapsto\Ind\epsilon_x$ is continuous.
Let
\[Y=\Ind^{-1}(U_2).\]
Then $Y$ is an open $G$-invariant neighbourhood of $z$ in $X$.  Note that  $x_n\in Y$ eventually.

Now suppose that, for some $y\in Y$, $G\cdot x_n\to G\cdot y$ in $Y/G$ and hence in $X/G$. Then  $\Ind\epsilon_{x_n}\to\Ind\epsilon_y$, and $\Ind\epsilon_y\in U_2$ since $y\in\Ind^{-1}(U_2)$.
But $(\Ind\epsilon_{x_n})$ has the unique limit $\Ind\epsilon_{z}$ in $U_2$, so $\Ind\epsilon_z=\Ind\epsilon_y$ and hence $\overline{G\cdot z} =\overline{G\cdot y}$ in $X$.  Since $G\cdot z$ is locally closed we obtain $G\cdot z=G\cdot y$ in $X$ and hence in $Y$.  (To see this, note that there exists a sequence $(r_k)_k\subset G$ such that $r_k\cdot y\to z$.  Let $U_3$ be an open subset of $X$ such that $G\cdot z=U_3\cap\overline{G\cdot z}$.  Then $r_k\cdot y\in U_3$ eventually, so eventually $r_k\cdot y\in U_3\cap\overline{G\cdot y}=U_3\cap\overline{G\cdot z}=G\cdot z$.)

Finally, we note that $C_0(Y)\times G$ is isomorphic to a closed two-sided
ideal $J$ of $C_0(X)\times G$ and $M_L(\Ind\epsilon_z,(\Ind\epsilon_{x_n}))$ is
the same whether we compute it in the ideal $J$ or in
$C_0(X)\times G$ (see \cite[Proposition 5.3]{ASS}); moreover,
$\phi_z^{-1}(V)= \phi_z^{-1}(V\cap Y)$ and $\phi_{x_n}^{-1}(V)=
\phi_{x_n}^{-1}(V\cap Y)$ when $n$ is large enough so that $x_n\in
Y$.  Thus we may replace $X$ by $Y$ and therefore assume that $G\cdot z$
is the unique limit of $G\cdot x_n$ in $X/G$, as claimed.

The idea of the rest of the proof is the same as in
Theorem~\ref{thm-Msquared}, but our estimates are more delicate. Fix
$\epsilon>0$ such that $M(1+\epsilon)^2<\lfloor M\rfloor +1$ and
choose $\gamma>0$ such that
\begin{equation}\label{eq-Mgamma}
\gamma<\frac{\epsilon\nu(\phi_z^{-1}(V))}{1+\epsilon}<\nu(\phi_z^{-1}(V)).
\end{equation}
By Lemma~\ref{lem-badV} there exists an open  relatively compact
neighbourhood $V_1$ of $z$ such that $\overline{V_1}\subset V$ and
\begin{equation*}
0<\nu(\phi_z^{-1}(V))-\gamma
<\nu(\phi_z^{-1}(V_1))
\leq\nu(\phi_z^{-1}(\overline{V_1}))
\leq\nu(\phi_z^{-1}(V))
<\nu(\phi_z^{-1}(V_1))+\gamma.
\end{equation*}
(The reason for passing from $V$ to $V_1$ is that we will  later
apply Lemma~\ref{lem-excision} to the compact neighbourhood
$\overline{V_1}$ and, in contrast to what could happen with
$\overline{V}$, we can control $\nu(\phi_z^{-1}(\overline{V_1}))$
relative to $\nu(\phi_z^{-1}(V_1))$.) Choose a subsequence
$(x_{n_i})_i$ of $(x_n)_n$ such that
\begin{equation*}
\nu(\phi_{x_{n_i}}^{-1}(V))\leq M\nu(\phi_z^{-1}(V))
\end{equation*}
for all $i\geq 1$.
Then
\begin{align}
\nu(\phi_{x_{n_i}}^{-1}(V_1))
&\leq\nu(\phi_{x_{n_i}}^{-1}(V))\notag\\
&\leq M\nu(\phi_z^{-1}(V))\notag\quad\quad\text{(by assumption)}\\
&<M\big( \nu(\phi_z^{-1}(V_1)) +\gamma \big)\notag\\
&<M\nu(\phi_z^{-1}(V_1)) +M\epsilon\big( \nu(\phi_z^{-1}(V))-\gamma \big)\notag\quad\quad\text{(by \eqref{eq-Mgamma})}\\
&<M\nu(\phi_z^{-1}(V_1)) +M\epsilon\nu(\phi_z^{-1}(V_1))\notag\\
&=M(1+\epsilon)\nu(\phi_z^{-1}(V_1))\label{eq-smallV}
\end{align}
for all $i$. Since
\[
\frac{\nu(\phi_z^{-1}(V_1))\big(\nu(\phi_z^{-1}(V_1))+\gamma+\frac{1}{j}  \big)}
{\big( \nu(\phi_z^{-1}(V_1))-\frac{1}{j} \big)^2}
\to 1+\frac{\gamma}{\nu(\phi_z^{-1}(V_1))}<1+\epsilon\]
as $j\to\infty$, there exists
$\delta>0$ such that $\delta<\nu(\phi_z^{-1}(V_1))$ and
\begin{equation}\label{eq-delta}
\frac{\nu(\phi_z^{-1}(V_1))\big(\nu(\phi_z^{-1}(\overline{V_1}))+\delta)
\big)}{\big( \nu(\phi_z^{-1}(V_1))-\delta \big)^2}
<\frac{\nu(\phi_z^{-1}(V_1))\big(\nu(\phi_z^{-1}(V_1))+\gamma+\delta)
\big)}{\big( \nu(\phi_z^{-1}(V_1))-\delta \big)^2}
<1+\epsilon.
\end{equation}

Next we construct a function $F\in C_c(X)$ with support inside
$V_1$. By the regularity of the measure $\nu$ there exists a
compact subset $W$ of the open set $\phi_z^{-1}(V_1)$ such that
$0<\nu(\phi_z^{-1}(V_1))-\delta<\nu(W)$. Since $W$ is compact, there
is a compact neighbourhood $W_1$ of $W$ contained in
$\phi_z^{-1}(V_1)$ and a continuous function $g:G\to[0,1]$ such that
$g$ is identically one on $W$ and is identically zero off the
interior of $W_1$. Then
\begin{equation}\label{eq-g}
\nu(\phi_z^{-1}(V_1))-\delta<\nu(W)\leq \int_G g(t)^2\, dt=\|g\|_2^2.
\end{equation}
There is a continuous function $g_1:W_1\cdot z\to[0,1]$ such that
$g_1(t\cdot z)=g(t)$ for $t\in W_1$. Since $W_1\cdot z$ is a
compact subset of $X$, we can extend $g_1$ to a continuous
function $g_2:X\to [0,1]$ using Tietze's Extension Theorem. There
exists a compact neighbourhood $P$ of $W_1\cdot z$ contained in
$V_1$ and a continuous function $h:X\to [0,1]$ such that $h$ is
identically one on $W_1\cdot z$ and is identically zero off the
interior of $P$.
 Note that $h$ has compact support contained in $P$. We set
$f(x)=h(x)g_2(x)$.
Then $f$ is continuous on $X$ and
\begin{equation}\label{eq-reciprocal2}
 \|f_z\|^2_2=\int_G f(t\cdot z)^2\, dt=
 \int_G h(t\cdot z)^2g_2(t\cdot z)^2\, dt
 \geq\int_{W_1} g(t)^2\, dt
 =\|g\|^2_2
\end{equation}
since $h$ is identically one on $W_1\cdot z$ and $g$ has  support
inside $W_1$. We now set $\displaystyle
{F(x)=\frac{f(x)}{\|f_z\|_2}} $. Thus $F\in C_c(X)$, $\|F_z\|_2=1$
and $F_x(s)=F(s\cdot x)\neq 0$ implies $s\in\phi_x^{-1}(V_1)$ by our
choice of $h$.

Let $K$ be an open relatively compact symmetric neighbourhood of
$(\supp F_z)(\supp F_z)^{-1}$ in $G$ and $L$ an open relatively
compact neighbourhood of $\supp F$ in $X$. Choose $b\in C_c(G\times
X)$ such that $0\leq b\leq 1$ and $b$ is identically one on the set
$(\supp F_z)(\supp F_z)^{-1}\times\supp F$ and $b$ is identically
zero off $K\times L$. (Thus $b$ is as in Theorem~\ref{thm-Msquared},
but we have rounded it off.) Set
\[
B(t,x)=F(x)F(t^{-1}\cdot x)b(t^{-1},x)\Delta(t)^{-1/2}\quad
\text{and}\quad D=\frac{1}{2}(B+B^*).
\]

Again,  $\Ind\epsilon_z(D)$, and hence $\Ind\epsilon_z(D^**D)$, is
the rank-one projection determined by the unit vector $F_z\in
L^2(G,\nu)$. From \eqref{eq-traceD} we have
\begin{align*}
\tr(&\Ind\epsilon_{x_{n_i}}(D^**D))\\
&=\frac{1}{4}\int_G F(s\cdot x_{n_i})^2\Big(  \int_G F(w\cdot
x_{n_i})^2\big( b(ws^{-1},s\cdot x_{n_i})+b(sw^{-1},w\cdot x_{n_i})
\big)^2\, d\nu(w)\Big)\, d\nu(s).
\end{align*}
The inner integrand is zero unless $w\in \phi_{x_{n_i}}^{-1}(V_1)\cap Ks$ because $F_{x_{n_i}}(w)=F(w\cdot x_{n_i})\neq 0$ implies $w\in\phi_{x_{n_i}}^{-1}(V_1)$ and $b$ is identically zero off $K\times L$. Thus
\begin{align*}
\tr(\Ind\epsilon_{x_{n_i}}(D^**D))
&\leq \int_{s\in\phi_{x_{n_i}}^{-1}(V_1)}F(s\cdot x_{n_i})^2\Bigg( \int_{w\in\phi_{x_{n_i}}^{-1}(V_1)\cap Ks} F(w\cdot x_{n_i})^2\, d\nu(w)\Bigg)\, d\nu(s)\\
&\leq \frac{1}{\|f_z\|_2^4}\int_{s\in\phi_{x_{n_i}}^{-1}(V_1)}1\Bigg( \int_{w\in\phi_{x_{n_i}}^{-1}(V_1)\cap Ks}1\, d\nu(w)\Bigg)\, d\nu(s).
\end{align*}

Choose an open neighbourhood $U$ of $\phi_z^{-1}(\overline{V_1})$ such that
$\nu(U)<\nu(\phi_z^{-1}(\overline{V_1}))+\delta$. By
Lemma~\ref{lem-excision}, applied with $\overline{V_1}$, $\overline{K}$ and $U$,
there exists $i_0$ such that, for every $i\geq
i_0$ and every $s\in \phi_{x_{n_i}}^{-1}(\overline{V_1})$ there exists
$r\in\phi_z^{-1}(\overline{V_1})$ with
$\overline{K}s\cap\phi_{x_{n_i}}^{-1}(\overline{V_1})\subset Ur^{-1}s$.
It follows that
\[
\nu\big( \overline{K}s\cap\phi_{x_{n_i}}^{-1}( \overline{V_1})
\big)\leq\nu(U)<\nu(\phi_z^{-1}(\overline{V_1}))+\delta
\]
by the right-invariance of $\nu$.
So, provided $i\geq i_0$,
\begin{align*}
\tr(\Ind\epsilon_{x_{n_i}}(D^**D))&\leq \frac{\nu(\phi_{x_{n_i}}^{-1}(V_1))\big(\nu(\phi_z^{-1}(\overline{V_1}))+\delta \big)}{\|f_z\|_2^4}\\
&<\frac{M(1+\epsilon)\nu(\phi_z^{-1}(V_1))\big(\nu(\phi_z^{-1}(\overline{V_1}))+\delta
\big)}{\|g\|_2^4}
\text{\quad\quad using \eqref{eq-smallV} and \eqref{eq-reciprocal2}}\\
&\leq\frac{M(1+\epsilon)\nu(\phi_z^{-1}(V_1))\big(\nu(\phi_z^{-1}(\overline{V_1}))+\delta
\big)}{(\nu(\phi_z^{-1}(V_1))-\delta)^2}
\text{\quad\quad using \eqref{eq-g}}\\
&<M(1+\epsilon)^2\text{\qquad\qquad\qquad\qquad\qquad\qquad\qquad using \eqref{eq-delta}}.
\end{align*}
By generalised lower semi-continuity \cite[Theorem~4.3]{AS}
\begin{align*}
\liminf_n \tr(\Ind\epsilon_{x_n}(D^**D))&\geq M_L(\Ind\epsilon_z,(\Ind\epsilon_{x_n}))\tr(\Ind\epsilon_z(D^**D))
\\&=M_L(\Ind\epsilon_z,(\Ind\epsilon_{x_n})).
\end{align*}
We now have
\begin{equation*}
M_L(\Ind\epsilon_z,(\Ind\epsilon_{x_n}))\leq \liminf_n \tr(\Ind\epsilon_{x_n}(D^**D))\leq M(1+\epsilon)^2<\lfloor M\rfloor+1,
\end{equation*} and hence $M_L(\Ind\epsilon_z,(\Ind\epsilon_{x_n}))\leq\lfloor M\rfloor$.
\end{proof}

The first corollary of Theorem~\ref{thm-M} is the analogous result for the upper multiplicity relative to a net.

\begin{thm}\label{thm-Mupper}
Suppose that $(G,X)$ is a free second-countable transformation
group. Let $M\in\R$ with $M\geq 1$, let $z\in X$  and let
$(x_n)_{n\geq 1}$ be a sequence in $X$. Assume that $G\cdot z$ is
locally closed in $X$. Suppose that there exists an open
neighbourhood $V$ of $z$ in $X$ such that $\phi_z^{-1}(V)$ is
relatively compact and
\begin{equation*}
\nu(\phi_{x_n}^{-1}(V))\leq M\nu(\phi_z^{-1}(V))<\infty
\end{equation*}
eventually. Then $M_U(\Ind\epsilon_z,(\Ind\epsilon_{x_n}))\leq\lfloor M\rfloor$.
\end{thm}

\begin{proof}
Both $G$ and $X$ are second countable, so $C_0(X)\times G$ is
separable. By Lemma~\ref{lemma-sepmult} there exists a
subsequence $(\Ind\epsilon_{x_{n_i}})_i$ such that
\[M_U(\Ind\epsilon_z,(\Ind\epsilon_{x_n}))=M_U(\Ind\epsilon_z,(\Ind\epsilon_{x_{n_i}}))=M_L(\Ind\epsilon_z,(\Ind\epsilon_{x_{n_i}})).\]
By Theorem~\ref{thm-M}, $M_L(\Ind\epsilon_z,(\Ind\epsilon_{x_{n_i}}))\leq \lfloor
M\rfloor$, and hence $M_U(\Ind\epsilon_z,(\Ind\epsilon_{x_n}))\leq \lfloor M\rfloor$.
\end{proof}

\begin{cor}\label{cor1-M}
Suppose that $(G,X)$ is a free second-countable transformation group
such that all the orbits are locally closed in $X$. Let $M\in\R$
with $M\geq 1$ and let $z\in X$. If for every sequence $(x_n)_{n\geq
1}$ in $X$ which converges to $z$ there exists an open neighbourhood
$V$ of $z$ in $X$ such that $\phi_z^{-1}(V)$ is relatively compact
and
\begin{equation*}
\nu(\phi_{x_n}^{-1}(V))\leq M\nu(\phi_z^{-1}(V))<\infty
\end{equation*}
frequently, then $M_U(\Ind\epsilon_z)\leq\lfloor M\rfloor$.
\end{cor}

\begin{proof}
Since  $C_0(X)\times G$ is separable  it follows from
\cite[Lemma~1.2]{AK} that there exists a sequence $(\pi_n)_{n\geq
1}$ in $(C_0(X)\times G)^\wedge$ converging to $\Ind\epsilon_z$, such that
\begin{equation*}
M_L(\Ind\epsilon_z,(\pi_n))=M_U(\Ind\epsilon_z,(\pi_n))=M_U(\Ind\epsilon_z).
\end{equation*}
Since the orbits are locally closed, $C_0(X)\times G$ is Type
\textrm{I} by \cite[Theorem~3.3]{Goot}, and then the map $x\mapsto
\Ind\epsilon_x$ induces a homeomorphism of $X/G$ onto $C_0(X)\times
G)^\wedge$ by
\cite[Lemma~16]{Green2}. In particular, the mapping
$X\to(C_0(X)\times G)^\wedge$ is an open surjection so there exists
a sequence $(x_i)_{i\geq 1}$ in $X$ converging to $z$ such that
$(\Ind\epsilon_{x_i})_{i\geq 1}$ is a subsequence of $(\pi_n)_{n\geq
1}$. By Theorem~\ref{thm-M},
$M_L(\Ind\epsilon_z,(\Ind\epsilon_{x_i}))\leq \lfloor M\rfloor$.
Since
\begin{align*}
M_U(\Ind\epsilon_z)&=M_L(\Ind\epsilon_z,(\pi_n))\leq M_L(\Ind\epsilon_z,(\Ind\epsilon_{x_i}))\leq M_U(\Ind\epsilon_z,(\Ind\epsilon_{x_i}))\\
&\leq M_U(\Ind\epsilon_z,(\pi_n))=M_U(\Ind\epsilon_z),
\end{align*}
we obtain $M_U(\Ind\epsilon_z)\leq\lfloor M\rfloor$.
\end{proof}

\section{Topological strength of convergence}\label{sec-top}


A sequence $(t_n)_{n\geq 1}\subset
G$ tends to infinity if  it admits no convergent subsequence. Let $k\in\P$.
 A sequence $(x_n)_{n\geq 1}$ in $X$ is
\emph{$k$-times convergent in $X/G$ to $z\in X$} if there exist $k$ sequences
$(t_n^{(1)})_n,(t_n^{(2)})_n,\cdots ,(t_n^{(k)})_n\subset G$, such
that
\begin{enumerate}
\item $t_n^{(i)}\cdot x_n\to z$ as $n\to\infty$ for $1\leq i\leq k$, and
\item if $1\leq i<j\leq k$ then $t_n^{(j)}(t_n^{(i)})^{-1}\to\infty$ as
$n\to\infty$.
\end{enumerate}
This definition of $k$-times convergence is a special case of \cite[Definition~2.2]{AD} obtained by taking $A=C_0(X)$ and replacing $t_n^{(i)}$ by $(t_n^{(i)})^{-1}$.

Proposition~\ref{prop-excision-combined} below shows that measure-theoretic accumulation implies topological strength of convergence in $X/G$; the proof uses our key technical Lemma~\ref{lem-excision}.

\begin{prop}\label{prop-excision-combined}
Let $(G, X)$ be a second-countable transformation group.  Let
$k\in\P$ and $z\in X$, with $G\cdot z$ locally closed in $X$.
Assume that $(x_n)_{n\geq 1}$ is a sequence in $X$ such that
$G\cdot x_n\to G\cdot z$ and $G\cdot z$ is the unique limit  of
$(G\cdot x_n)_n$ in $X/G$.
\begin{enumerate}
\item Suppose that there
exists a basic sequence $(W_m)_{m\geq 1}$ of compact neighbourhoods of $z$
with $W_{m+1}\subset W_m$, such that, for each $m$
\[
\liminf_n\nu(\phi_{x_n}^{-1}(W_m))>(k-1)\nu(\phi_z^{-1}(W_m)).
\]
Then  $(x_n)$  converges $k$-times in $X/G$ to $z$.
\item Suppose that there
exists a basic sequence $(W_m)_{m\geq 1}$ of compact neighbourhoods of $z$
with $W_{m+1}\subset W_m$, such that, for each $m$
\begin{equation*}
\limsup_n\nu(\phi_{x_n}^{-1}(W_m))>(k-1)\nu(\phi_z^{-1}(W_m)).
\end{equation*}
Then  there exists a subsequence of $(x_n)_n$ which  converges $k$-times in $X/G$ to $z$ .
\end{enumerate}
\end{prop}

\begin{proof}
Let $(K_m)_{m\geq 1}$ be an increasing sequence of  compact subsets of $G$ such
that $G=\cup_{m\geq 1}\interior(K_m)$.

(1) By regularity of $\nu$, for each $m\geq 1$, there exists an
open neighbourhood $U_m$ of $\phi_z^{-1}(W_m)$ such that
\begin{equation}\label{eq-reg}
\liminf_n\nu(\phi_{x_n}^{-1}(W_m))>(k-1)\nu(U_m).
\end{equation}
We will  construct, by induction,  a strictly increasing sequence of
positive integers $(n_m)_{m\geq 1}$ such that, for all $n\geq
n_m$,
\begin{align}
\nu(K_ms\cap\phi_{x_n}^{-1}(W_m))&\leq\nu(U_m)\quad\text{for all\
}s\in\phi_{x_n}^{-1}(W_m),
\text{\ and}\label{eq-excised}\\
\nu(\phi_{x_n}^{-1}(W_m))&>(k-1)\nu(U_m).
\label{eq-measure}
\end{align}
We construct $n_1$ by applying Lemma~\ref{lem-excision} to
$K_1, W_1$ and $U_1$ to obtain \eqref{eq-excised} for $m=1$, and then,
if necessary, increasing $n_1$ to ensure that \eqref{eq-measure}
holds using \eqref{eq-reg} with $m=1$.
Assuming that we have constructed $n_1<n_2<\dots<n_{m-1}$,
we apply Lemma~\ref{lem-excision} to $K_m, W_m$ and $U_m$
to obtain $n_m>n_{m-1}$ such that \eqref{eq-excised} holds,
and again, if necessary,  increase $n_m$ to obtain \eqref{eq-measure}.

If $n_1>1$ then, for $1\leq n<n_1$, we set $t^{(i)}_n=e$ for $1\leq i\leq k$.
For each $n\geq n_1$ there is a unique $m$ such that $n_m\leq n< n_{m+1}$.
Choose $t^{(1)}_n\in\phi_{x_n}^{-1}(W_m)$.  Using \eqref{eq-excised} and \eqref{eq-measure}
\[
\nu(\phi_{x_n}^{-1}(W_m)\setminus K_mt^{(1)}_n)>(k-2)\nu(U_m).
\]
So we may choose $t_n^{(2)}\in \phi_{x_n}^{-1}(W_m)\setminus K_mt^{(1)}_n$.
Continuing in this way, we use \eqref{eq-excised} and \eqref{eq-measure} to choose
\begin{align*}
t_n^{(3)}&\in \phi_{x_n}^{-1}(W_m)\setminus (K_mt^{(1)}_n\cup K_mt^{(2)}_n)\\
\vdots & \\
t_n^{(k)}&\in \phi_{x_n}^{-1}(W_m)\setminus (\cup_{i=1}^{k-1}K_mt^{(i)}_n).
\end{align*}
Note that for $n_m\leq n<n_{m+1}$ we have
\[
t_n^{(i)}\cdot x_n\in W_m\text{\ for\ }1\leq i\leq k
\quad\text{and}\quad
t_n^{(j)}\notin K_mt_n^{(i)}\text{\ for\ }1\leq i<j\leq k.
\]

We claim that $t_n^{(i)}\cdot x_n\to z$ as $n\to\infty$ for $1\leq i\leq k$.
To see this, fix $i$ and let $V$ be a neighbourhood of $z$.
There exists $m_0$ such that $W_m\subset V$ for all $m\geq m_0$.
For each $n\geq n_{m_0}$ there exists a  $m\geq m_0$ such that $n_m\leq n<n_{m+1}$,
and thus $t_n^{(i)}\cdot x_n\in W_m\subset V$.

Next,  we claim that $t_n^{(j)}(t_n^{(i)})^{-1}\to\infty$ as $n\to\infty$ for $1\leq
i<j\leq k$. Fix $i<j$ and let $K$ be a compact
subset of $G$. There exists $m_0$ such that $K_m\supset K$ for all
$m\geq m_0$. Then for $n\geq n_{m_0}$ there exists a $m\geq m_0$
such that $n_m\leq n<n_{m+1}$, and hence
$t_n^{(j)}(t_n^{(i)})^{-1}\in G\setminus K_m\subset G\setminus K$.
We have shown that $(x_n)_n$ converges $k$-times in $X/G$  to $z$.

\medskip
(2) By regularity of
$\nu$, for each $m\geq 1$ there exists an open neighbourhood $U_m$ of
$\phi_z^{-1}(W_m)$ such that
\begin{equation}\label{eq-reg2}
\limsup_n\nu(\phi_{x_n}^{-1}(W_m))>(k-1)\nu(U_m).
\end{equation}
We will  construct, by induction,  an increasing sequence $(i_m)_{m\geq
1}$ such that
\begin{gather}
\nu(K_{m}s\cap\phi_{x_{i_m}}^{-1}(W_{m}))\leq\nu(U_{m})\text{\
for all\
}s\in\phi_{x_{i_m}}^{-1}(W_{m})\text{\ and}\label{eq-excisedm}\\
\nu(\phi_{x_{i_m}}^{-1}(W_m))>(k-1)\nu(U_m)\label{eq-moveup}.
\end{gather}
To start, we apply Lemma~\ref{lem-excision} to $K_1, W_1$ and $U_1$
to obtain $n_1$ such that
\begin{equation*}
n\geq n_1\text{\ implies\
}\nu(K_1s\cap\phi_{x_n}^{-1}(W_1))\leq\nu(U_1)\text{\ for all\
}s\in\phi_{x_n}^{-1}(W_1).
\end{equation*}
Then, using \eqref{eq-reg2} with $m=1$, choose $i_1\geq n_1$ such
that
\begin{equation*}
\nu(\phi_{x_{i_1}}^{-1}(W_1))>(k-1)\nu(U_1).
\end{equation*}
Assuming that $i_1<i_2<\dots<i_{m-1}$ have been chosen, apply
Lemma~\ref{lem-excision} to $K_m, W_m$ and $U_m$ to obtain
$n_m>i_{m-1}$ such that
\begin{equation*}
n\geq n_m\text{\ implies\
}\nu(K_ms\cap\phi_{x_n}^{-1}(W_m))\leq\nu(U_m)\text{\ for all\
}s\in\phi_{x_n}^{-1}(W_m),
\end{equation*}
and then, using \eqref{eq-reg2}, choose $i_m\geq n_m$ such that
\begin{equation*}
\nu(\phi_{x_{i_m}}^{-1}(W_m))>(k-1)\nu(U_m).
\end{equation*}

Now, for each $m\geq 1$, choose
$t^{(1)}_{i_m}\in\phi_{x_{i_m}}^{-1}(W_m)$. Using
\eqref{eq-excisedm} and \eqref{eq-moveup}
\[
\nu(\phi_{x_{i_m}}^{-1}(W_m)\setminus
K_mt^{(1)}_{i_m})>(k-2)\nu(U_m).
\]
So we may choose $t_{i_m}^{(2)}\in
\phi_{x_{i_m}}^{-1}(W_{m})\setminus K_{m}t^{(1)}_{i_m}$. Continuing
in this way, we use \eqref{eq-excisedm} and \eqref{eq-moveup} to
choose
\begin{align*}
t_{i_m}^{(3)}&\in \phi_{x_{i_m}}^{-1}(W_m)\setminus  (K_mt^{(1)}_{i_m}\cup K_mt^{(2)}_{i_m})\\
\vdots & \\
t_{i_m}^{(k)}&\in \phi_{x_{i_m}}^{-1}(W_m)\setminus
(\cup_{j=1}^{k-1}K_mt^{(j)}_{i_m}).
\end{align*}
Since $t^{(j)}_{i_m}\cdot x_{i_m}\in W_m$ and $(W_m)_m$ is
decreasing, it follows that $t^{(j)}_{i_m}\cdot x_{i_m}\to_m z$ for
$1\leq j\leq k$.  If $K$ is any compact subset of $G$ then there
exists $m_0$ such that $K\subset K_m$ for all $m\geq m_0$. Then for
$m\geq m_0$ we have $t^{(j)}_{i_m}(t^{(l)}_{i_m})^{-1}\in G\setminus
K_m\subset G\setminus K$ for $1\leq l<j \leq k$. Thus $(x_{i_m})_m$ is a
subsequence of $(x_n)_n$ which converges $k$-times  in $X/G$ to $z$.
\end{proof}


Although we shall improve Proposition~\ref{propcor-msq} later (see
Corollary~\ref{cor1-main}), we prove it now because it will be
needed in the proof of Theorem~\ref{thm-main}.


\begin{prop}\label{propcor-msq}
Suppose that $(G,X)$ is a free second-countable transformation
group. Let $z\in X$ and let $(x_n)_{n\geq1}$ be a sequence in $X$.
Assume that $G\cdot z$ is locally closed in $X$.  Consider the
following properties.
\begin{enumerate}
\item
$M_L(\Ind\epsilon_z,(\Ind\epsilon_{x_n}))=\infty$;
 \item for every open neighbourhood $V$ of $z$ such that
$\phi^{-1}_z(V)$ is relatively compact, $\nu(\phi^{-1}_{x_n}(V))\to
\infty$ as $n\to\infty$;
\item for each $k\geq1$, the sequence $(x_n)_n$
converges $k$-times  in $X/G$ to $z$.
\end{enumerate}
Then (1) $\Longrightarrow$ (2) and (2) $\Longrightarrow$ (3).
\end{prop}

\begin{proof}
(1) $\Longrightarrow$ (2). This is an immediate corollary of
Theorem~\ref{thm-Msquared}.

(2)$\Longrightarrow$ (3). Let $(K_m)_{m\geq 1}$ be an increasing
sequence of compact subsets of $G$ such that $G=\cup_{m\geq
1}\interior(K_m)$ and let $(V_m)_{m\geq1}$ be a decreasing sequence
of open,  basic neighbourhoods of $z$ such  that $\phi_z^{-1}(V_1)$
is relatively compact (such neighbourhoods exist by
Lemma~\ref{lem-remark}). Let $k\geq1$. Assuming (2), we have
$$\nu(\phi^{-1}_{x_n}(V_m)) \to_n\infty \qquad (m\geq1).$$
Hence we can construct inductively a strictly increasing sequence
of positive integers $(n_m)_{m\geq 1}$ such that, for all $n\geq
n_m$,
\begin{equation}\label{big-measure}
\nu(\phi_{x_n}^{-1}(V_m))>(k-1)\nu(K_m).
\end{equation}
We can now choose $t^{(i)}_n$ for $1\leq i\leq k$ as
in Proposition~\ref{prop-excision-combined}(1) by using
\eqref{big-measure} in place of  \eqref{eq-excised} and
\eqref{eq-measure}.
\end{proof}


\section{Main results}\label{sec-main}

We will shortly combine the results from  \S\ref{sec-measure} and
\S\ref{sec-top} to prove our main theorem stated in the
introduction.

\begin{lemma}\label{lem-C} Suppose that $(G,X)$ is a transformation
group.  Let $k\in \P$, $z\in X$ and $(x_n)_{n\geq 1}$ be a
sequence in $X$. Assume that $G\cdot z$ is locally closed in $X$
and that there exists a real number $R>k-1$ such that for every
open neighbourhood $U$ of $z$ with $\phi_z^{-1}(U)$ relatively
compact we have
\[
\liminf_n \nu(\phi_{x_n}^{-1}(U))\geq R\nu(\phi_z^{-1}(U)).
\]
Given an open  neighbourhood $V$ of $z$ such that $\phi_z^{-1}(V)$
is relatively compact, there exists a compact neighbourhood $N$ of
$z$ with $N\subset V$ such that
\[
\liminf_n \nu(\phi_{x_n}^{-1}(N))> (k-1)\nu(\phi_z^{-1}(N)).
\]
\end{lemma}

\begin{proof}
Apply Lemma~\ref{lem-badV} to $V$ with
$0<\gamma<\big(\frac{R-k+1}{R}\big)\nu(\phi_z^{-1}(V))$ to get an open
relatively compact neighbourhood $V_1$ of $z$ with
$\overline{V_1}\subset V$ and
\begin{equation*}
\nu(\phi_z^{-1}(V))-\gamma <\nu(\phi_z^{-1}(V_1))
\leq\nu(\phi_z^{-1}(\overline{V}_1)) \leq\nu(\phi_z^{-1}(V))
<\nu(\phi_z^{-1}(V_1))+\gamma.
\end{equation*}
Since $\phi_z^{-1}(V_1)$ is relatively compact we have
\begin{align*}
\liminf_n\nu(\phi_{x_n}^{-1}(\overline{V}_1))
&\geq \liminf_n\nu(\phi_{x_n}^{-1}(V_1))\\
&\geq R\nu(\phi_z^{-1}(V_1)) \text{\quad by hypothesis}\\
&>R\big(\nu(\phi_z^{-1}(V))-\gamma \big)\\
&>(k-1)\nu(\phi_z^{-1}(V))\text{\quad  by our choice of $\gamma$}\\
&\geq (k-1)\nu(\phi_z^{-1}(\overline{V}_1)).
\end{align*}
So we may take $N=\overline{V}_1$.
\end{proof}

\begin{remark}\label{rem-lem-C}
An examination of the proof of Lemma~\ref{lem-C} shows that there is
a variant of the lemma in which $\liminf$ is replaced by $\limsup$.
\end{remark}

\begin{proof}[Proof of Theorem~\ref{thm-main}]
(1) $\Longrightarrow$ (2). Let $\phi$ be a pure state associated
with $\Ind\epsilon_z$ and  $\N$  the usual weak*-neighbourhood base
at zero in $(C_0(X)\times G)^*$. Suppose that $M_L(\Ind\epsilon_z,
(\Ind\epsilon_{x_n}))=r<k$.  There exists $N\in\N$ such that
\[
\liminf_n d((\Ind\epsilon_{x_n}),\phi, N)=r
\]
and hence a subsequence $(x_{n_i})_i$ such that
\begin{equation}\label{eq-mlessk}
d((\Ind\epsilon_{x_{n_i}}),\phi, N)=r<k
\end{equation}
for all $i\geq 1$. Note that $(x_{n_i})_i$ converges $k$-times  in
$X/G$ to $z$ because $(x_n)_n$ does. The proof of
\cite[Theorem~2.3]{AD}\footnote{The hypothesis in \cite[Theorem~2.3]{AD} that all the orbits
are locally closed in $\hat A$  is not
needed in the case $A=C_0(X)$, since $\Ind\epsilon_x$ is irreducible
by \cite[Proposition~4.2]{W}. The representation
$\tilde\epsilon_{x_{n_i}}\times\lambda$ on $L^2(G,\mu)$ used in
\cite[Theorem~2.3]{AD} is unitarily equivalent to our
$\Ind\epsilon_{x_{n_i}}$ on $L^2(G,\nu)$ via the unitary
$W:L^2(G,\mu)\to L^2(G,\nu)$ defined by
$(W\xi)(s)=\Delta(s)^{1/2}\xi(s)$.} shows that, after passing to a
further subsequence and reindexing, there exist unit vectors
$\eta_i^{(j)}$ $(1\leq j\leq k)$ in $L^2(G, \nu)$ such that
\[
\lim_{i\to\infty}\langle\Ind\epsilon_{x_{n_i}}(\cdot)\eta_i^{(j)}\, ,\, \eta_i^{(j)}\rangle=\phi\quad(1\leq j\leq k),
\]
and that there exists $i_0$ such that if $i\geq i_0$ then
$\eta_i^{(j)}$ and $\eta_i^{(l)}$ are orthogonal if $j\neq l$.  By
increasing $i_0$ if necessary, we may assume that
$\langle\Ind\epsilon_{x_{n_i}}(\cdot)\eta_i^{(j)}\, ,\,
\eta_i^{(j)}\rangle\in\phi+N$ for all $i\geq i_0$ and $1\leq j\leq
k$. Thus $d((\Ind\epsilon_{x_{n_i}}),\phi, N)\geq k$ for all $i\geq
i_0$, contradicting \eqref{eq-mlessk}.

(2) $\Longrightarrow$ (3). If $M_L(\Ind\epsilon_z,(\Ind\epsilon_{x_n}))\geq k$ then
$M_L(\Ind\epsilon_z,(\Ind\epsilon_{x_n}))>\lfloor k-\epsilon\rfloor$ for every
$\epsilon>0$.  By Theorem~\ref{thm-M}, for every open
neighbourhood $V$ of $z$ such that
 $\phi_z^{-1}(V)$ is relatively compact,
\[
\nu(\phi_{x_n}^{-1}(V))> (k-\epsilon)\nu(\phi_z^{-1}(V))
\]
eventually, and hence (3) holds.

(3)   $\Longrightarrow$ (4) is immediate.

(4) $\Longrightarrow$ (5). Let $(V_j)_{j\geq 1}$ be a decreasing
sequence of basic  open neighbourhoods of $z$ such that
$\phi_z^{-1}(V_1)$ is relatively compact (such neighbourhoods exist
by Lemma~\ref{lem-remark}). By Lemma~\ref{lem-C} there exists a
compact neighbourhood $W_1$ of $z$ such that $W_1\subset V_1$ and
\[
\nu(\phi_{x_n}^{-1}(W_1))> (k-1)\nu(\phi_z^{-1}(W_1)).
\]
Now assume there are compact neighbourhoods  $W_1, W_2, \dots, W_m$
of $z$ with $W_1\supset W_2\supset\cdots \supset W_m$  such that
\begin{equation}\label{eq-3gives4}
W_i\subset V_i\text{\ and\ }\nu(\phi_{x_n}^{-1}(W_i))>
(k-1)\nu(\phi_z^{-1}(W_i))
\end{equation}
for $1\leq i\leq m$. Apply Lemma~\ref{lem-C} to  $(\interior
W_m)\cap V_{m+1}$ to obtain a compact neighbourhood $W_{m+1}$ of $z$
such that $W_{m+1}\subset (\interior W_m)\cap V_{m+1}$ and
\eqref{eq-3gives4} holds for $i=m+1$.

(5) $\Longrightarrow$ (1). We show first that $G\cdot x_n\to G\cdot
z$ in $X/G$.  Let $q:X\to X/G$ be the quotient map. Let $U$ be a
neighbourhood of $G\cdot z$ in $X/G$ and $V=q^{-1}(U)$. There exists
$m$ such that $W_m\subset V$. Since
$\liminf_n\nu(\phi_{x_n}^{-1}(W_m))>0$ there exists $n_0$ such that
$\phi_{x_n}^{-1}(W_m)\neq\emptyset$ for $n\geq n_0$. Thus, for
$n\geq n_0$,
\[G\cdot x_n=q(x_n)\in q(W_m)\subset q(V)= U.\]
Thus $G\cdot x_n$ is eventually in every neighbourhood  of $G\cdot
z$ in $X/G$.

Now suppose that
$M_L(\Ind\epsilon_z,(\Ind\epsilon_{x_n}))<\infty$. Then, as in the proof of
Theorem~\ref{thm-M}, we may localise to an open $G$-invariant
neighbourhood $Y$ of $z$ such that $G\cdot z$ is the unique limit
in $Y/G$ of the sequence $(G\cdot x_n)_n$.  Eventually $W_m\subset
Y$, and so the sequence $(x_n)_n$ converges $k$-times in
$Y/G$  to $z$ by Proposition~\ref{prop-excision-combined}(1)
applied to $Y$. But now $(x_n)_n$ converges $k$-times  in $X/G$ to $z$ as well.

Finally, if
$M_L(\Ind\epsilon_z,(\Ind\epsilon_{x_n}))=\infty$ then $(x_n)_n$ converges $k$-times
 in $X/G$ to $z$ by Proposition~\ref{propcor-msq}.
\end{proof}

\begin{thm}\label{thm-main-upper}
Suppose that $(G,X)$ is a free second-countable  transformation
group. Let $k\in \P$, let $z\in X$ and  let $(x_n)_{n\geq 1}\subset
X$ be a sequence such that $(G\cdot x_n)_n$ converges to $G\cdot z$
in $X/G$. Assume that $G\cdot z$ is locally closed in $X$. Then the
following are equivalent:
\begin{enumerate}
\item there exists a subsequence  $(x_{n_i})_{i\geq 1}$
of $(x_n)$ which converges $k$-times  in $X/G$ to $z$;
\item $M_U(\Ind\epsilon_z,(\Ind\epsilon_{x_n}))\geq k$;
\item for every open  neighbourhood $V$
of $z$ such that $\phi_z^{-1}(V)$ is relatively compact we have
\[
\limsup_n\nu(\phi_{x_n}^{-1}(V))\geq k\nu(\phi_z^{-1}(V));
\]
\item there exists a real number $R>k-1$ such that for every open
neighbourhood $V$ of $z$  with $\phi_z^{-1}(V)$ relatively compact
we have
\[
\limsup_n\nu(\phi_{x_n}^{-1}(V))\geq R\nu(\phi_z^{-1}(V));\]
\item there exists a decreasing sequence of basic compact
neighbourhoods $(W_m)_{m\geq 1}$ of $z$ such that, for each
$m\geq1$,
\[
\limsup_n\nu(\phi_{x_n}^{-1}(W_m))> (k-1)\nu(\phi_z^{-1}(W_m)).
\]
\end{enumerate}
\end{thm}

\begin{proof}
If  (1) holds then
$M_L(\Ind\epsilon_z,(\Ind\epsilon_{x_{n_i}}))\geq k$  using
Theorem~\ref{thm-main}, and hence
\[
M_U(\Ind\epsilon_z,(\Ind\epsilon_{x_n}))\geq
M_U(\Ind\epsilon_z,(\Ind\epsilon_{x_{n_i}}))\geq
M_L(\Ind\epsilon_z,(\Ind\epsilon_{x_{n_i}}))\geq k.
\]

If (2) holds then by Lemma~\ref{lemma-sepmult} there is a
subsequence $(x_{n_r})_r$ such that
$M_L(\Ind\epsilon_z,(\Ind\epsilon_{x_{n_r}}))=M_U(\Ind\epsilon_z,(\Ind\epsilon_{x_n}))$
so that $M_L(\Ind\epsilon_z,(\Ind\epsilon_{x_{n_r}}))\geq k$. Let
$V$ be any open neighbourhood of $z$ such that $\phi_z^{-1}(V)$ is
relatively compact. Then
\begin{equation*}
\limsup_n\nu(\phi_{x_n}^{-1}(V))\geq
\limsup_r\nu(\phi_{x_{n_r}}^{-1}(V)) \geq
\liminf_r\nu(\phi_{x_{n_r}}^{-1}(V)) \geq k\nu(\phi_{z}^{-1}(V)),
\end{equation*}
where we have used Theorem~\ref{thm-main} at the last step.

That (3) implies (4) is immediate.  That (4) implies (5) follows
from Remark~\ref{rem-lem-C} as in the proof of
Theorem~\ref{thm-main}.

Assume (5).  First suppose that
$M_L(\Ind\epsilon_z,(\Ind\epsilon_{x_n}))<\infty$. Since $G\cdot
x_n\to_n G\cdot z$,  we can localise to an open $G$-invariant
neighbourhood  $Y$ of $z$ such that $G\cdot z$ is the unique limit
of $(G\cdot x_n)_n$, as in the proof of Theorem~\ref{thm-M}. Now (1)
follows by applying Proposition~\ref{prop-excision-combined}(2) to
$Y$.

If $M_L(\Ind\epsilon_z,(\Ind\epsilon_{x_n}))=\infty$ then $(x_n)_n$
converges $k$-times
 in $X/G$ to $z$ by Proposition~\ref{propcor-msq}.
\end{proof}

We now derive some further consequences of Theorem~\ref{thm-main}


\begin{cor}\label{cor-converseAD}
Suppose that $(G,X)$ is a free second-countable  transformation
group such that all the orbits are locally closed in $X$. Let $k\in\P$
and let $z\in X$. Then the following are equivalent:
\begin{enumerate}
\item there exists a sequence $(x_n)_{n\geq 1}$ in $X$ which is
$k$-times convergent  in $X/G$ to $z$;
\item $M_U(\Ind\epsilon_z)\geq k$.
\end{enumerate}
\end{cor}

\begin{proof}
Assume (1). Then  $M_U(\Ind\epsilon_z)\geq M_L(\Ind\epsilon_z,(\Ind\epsilon_{x_n}))\geq k$ by Theorem~\ref{thm-main}.

Assume (2). By \cite[Lemma~1.2]{AK} there exists a sequence
$(\pi_n)_{n\geq 1}$ converging to $\Ind\epsilon_z$ such that
$M_L(\Ind\epsilon_z,(\pi_n))=M_U(\Ind\epsilon_z,(\pi_n))=M_U(\Ind\epsilon_z)$.
 Since the orbits are locally closed, the mapping $X\to(C_0(X)\times G)^\wedge:x\mapsto \Ind\epsilon_x$ is
a  surjection. So there is a sequence $(x_n)_n$ in $X$ such that
$\Ind\epsilon_{x_n}=\pi_n$ for each $n$. Then
\[M_L(\Ind\epsilon_z,(\Ind\epsilon_{x_n}))= M_L(\Ind\epsilon_z,(\pi_n))=M_U(\Ind\epsilon_z)\geq k,\]
and it follows from Theorem~\ref{thm-main} that $(x_n)_{n\geq 1}$ is
$k$-times convergent  in $X/G$ to $z$.
\end{proof}

That (1) $\Longrightarrow$ (2) in  Corollary~\ref{cor-converseAD}
 is a special case of \cite[Theorem~2.3]{AD}.

\begin{cor}\label{cor-3-main2}
Suppose that $(G,X)$ is a free second-countable transformation group such that all the
orbits are locally closed in $X$.  Let $k\in\P$ and let $z\in X$ such that
$G\cdot z$ is  not open in $X$. Then the following are
equivalent:
\begin{enumerate}
\item whenever $(x_n)_{n\geq 1}$ is a sequence in $X$ which converges to $z$ and
satisfies $z\notin G\cdot x_n$ eventually, then $(x_n)$ is $k$-times convergent  in $X/G$ to $z$;
\item $M_L(\Ind\epsilon_z)\geq k$.
\end{enumerate}
\end{cor}

\begin{proof}
Assume (1).  By Lemma~\ref{lem-Lsequence}, there is a sequence
$(\pi_n)_{n\geq 1}$ in $(C_0(X)\times G)^\wedge$ such that
$\pi_n\neq\Ind\epsilon_z$ for all $n$, $\pi_n\to_n\Ind\epsilon_z$ and
\begin{equation}\label{eq-lower}M_L(\Ind\epsilon_z)=M_L(\Ind\epsilon_z,(\pi_n)) =M_U(\Ind\epsilon_z,(\pi_n)).\end{equation}
Since the orbits are locally closed, the mapping
$X\to(C_0(X)\times G)^\wedge:x\mapsto \Ind\epsilon_x$ is an open
surjection. So there is a sequence $(x_n)_n$ in $X$ such that
$x_n\to z$ and $(\Ind\epsilon_{x_n})_n$ is a subsequence of
$(\pi_n)_n$. Using \eqref{eq-lower},
$M_L(\Ind\epsilon_z)=M_L(\Ind\epsilon_z,(\Ind\epsilon_{x_n}))$. By
(1), $(x_n)_n$ is $k$-times convergent  in $X/G$ to $z$, so it
follows from  Theorem~\ref{thm-main} that
$M_L(\Ind\epsilon_z)=M_L(\Ind\epsilon_z,(\Ind\epsilon_{x_n}))\geq
k$.

Assume (2). If $(x_n)_n$ is a sequence in $X$ which converges to $z$ such that $z\notin G\cdot x_n$ eventually, then
\[M_L(\Ind\epsilon_x,(\Ind\epsilon_{x_n}))\geq M_L(\Ind\epsilon_z)\geq k.\]
 It  follows from Theorem~\ref{thm-main} that $(x_n)_n$  is $k$-times convergent in $X/G$ to $z$.
\end{proof}

Corollary~\ref{cor1-main} improves Proposition~\ref{propcor-msq}, and
is immediate from Theorem~\ref{thm-main}.

\begin{cor}\label{cor1-main}
Suppose that $(G,X)$ is a free second-countable transformation
group. Let $z\in X$ and let $(x_n)_{n\geq1}$ be a sequence in $X$.
Assume that $G\cdot z$ is locally closed in $X$.  Then the following
are equivalent:
\begin{enumerate}
\item for each $k\geq1$, the sequence $(x_n)_n$ converges $k$-times
 in $X/G$ to $z$;
\item $M_L(\Ind\epsilon_z,(\Ind\epsilon_{x_n}))=\infty$;
\item for
every open neighbourhood $V$ of $z$ such that $\phi_z^{-1}(V)$ is
relatively compact, $\nu(\phi^{-1}_{x_n}(V))\to \infty$ as
$n\to\infty$.
\end{enumerate}
\end{cor}


\begin{cor}\label{cor2a-main-new1}
Suppose that $(G,X)$ is a free second-countable transformation
group. Let $z\in X$ and let $(x_n)_{n\geq1}\subset X$ be a sequence
converging to $z$.  Assume that $G\cdot z$ is locally closed in $X$.
Then the following are equivalent:
\begin{enumerate}
\item there exists an open
neighbourhood $V$ of $z$ such that $\phi^{-1}_z(V)$ is relatively
compact and \[\limsup_n\nu(\phi^{-1}_{x_n}(V))<\infty;\]
\item $M_U(\Ind\epsilon_z,(\Ind\epsilon_{x_n}))<\infty$.
\end{enumerate}
\end{cor}

\begin{proof}
Suppose that (1) holds. Since $C_0(X)\times G$ is separable, it
follows from Lemma~\ref{lemma-sepmult} that there exists a
subsequence $(x_{n_j})_{j\geq1}$ of $(x_n)$ such that
$$M_L(\Ind\epsilon_z,(\Ind\epsilon_{x_{n_j}})) = M_U(\Ind\epsilon_z,(\Ind\epsilon_{x_{n_j}})) =
M_U(\Ind\epsilon_z,(\Ind\epsilon_{x_n})).$$  By (1) and
Corollary~\ref{cor1-main},
$M_L(\Ind\epsilon_z,(\Ind\epsilon_{x_{n_j}}))<\infty$. Hence
$M_U(\Ind\epsilon_z,(\Ind\epsilon_{x_n}))<\infty$, that is (2)
holds.

Suppose that (1) fails. Let $(V_i)_{i\geq1}$ be a decreasing
sequence of open basic neighbourhoods of $z$ such that
$\phi_z^{-1}(V_1)$ is relatively compact (such neighbourhoods exist
by Lemma~\ref{lem-remark}). Then
$\limsup_n\{\nu(\phi^{-1}_{x_n}(V_i))\}=\infty$ for each $i$ and we
may choose a subsequence $(x_{n_i})_i$ of $(x_n)_n$ such that
$\nu(\phi^{-1}_{x_{n_i}}(V_i))\to_i\infty$.

Let $V$ be any open neighbourhood of $z$ such that $\phi^{-1}_z(V)$
is relatively compact. There exists $i_0$ such that $V_i\subset V$
for all $i\geq i_0$. Then, for $i\geq i_0$,
$$
\nu(\phi^{-1}_{x_{n_i}}(V_i)) \leq
\nu(\phi^{-1}_{x_{n_i}}(V)).
$$
Thus $\nu(\phi^{-1}_{x_{n_i}}(V))\to_i\infty$. By
Corollary~\ref{cor1-main},
$M_L(\Ind\epsilon_z,(\Ind\epsilon_{x_{n_i}}))=\infty$. Hence
$M_U(\Ind\epsilon_z,(\Ind\epsilon_{x_n}))=\infty$, that is (2)
fails.
\end{proof}

\begin{cor}\label{cor2a-main-new2}
Suppose that $(G,X)$ is a free second-countable transformation
group such that all the orbits are locally closed in $X$. Let $z\in X$
and let $(x_n)_{n\geq1}\subset X$ be a sequence converging to $z$.
Then the following are equivalent:
\begin{enumerate}
\item $M_U(\Ind\epsilon_z)<\infty$;
\item there exists an open
neighbourhood $V$ of $z$ such that $\phi^{-1}_z(V)$ is relatively
compact and \[ \sup_{x\in V}\nu(\phi^{-1}_x(V))<\infty;
\]
\end{enumerate}
\end{cor}
\begin{proof}
If (2) holds then (1) holds by Corollary~\ref{cor1-M}.

Let $(V_i)_{i\geq1}$ be a decreasing sequence of  open basic
neighbourhoods of $z$ such that $\phi_z^{-1}(V_1)$ is relatively
compact. If (2) fails then $\sup_{x\in
V_i}\{\nu(\phi^{-1}_x(V_i))\}=\infty$ for each $i$ and  we may
choose a sequence $(x_i)_i$ such that $x_i\in V_i$ for all $i$ and
$\nu(\phi^{-1}_{x_i}(V_i))\to_i\infty$. Since $(V_i)_i$ is
decreasing, $x_i\to_i z$.

Let $V$ be an open neighbourhood of $z$ such that $\phi^{-1}_z(V)$
is relatively compact. There exists $i_0$ such that $V_i\subset V$
for all $i\geq i_0$. Then, for $i\geq i_0$,
$$
\nu(\phi^{-1}_{x_i}(V_i)) \leq \nu(\phi^{-1}_{x_i}(V)).
$$
Thus $\nu(\phi^{-1}_{x_i}(V))\to_i\infty$. By  the contrapositive of
the (2) $\Longrightarrow$ (1) direction of
Corollary~\ref{cor2a-main-new1},
$M_U(\Ind\epsilon_z,(\Ind\epsilon_{x_i}))=\infty$. Hence
$M_U(\Ind\epsilon_z)=\infty$, that is (1) fails.
\end{proof}

In the terminology of \cite{aH2},  property (2) of
Corollary~\ref{cor2a-main-new2} says that $V$ is an integrable
neighbourhood of $z$; in view of \cite[Theorem~2.8]{ASS}, the
equivalence of (1) and (2) is essentially contained in
\cite[Theorem~5.8(2)]{aH2}.

\begin{remark}
 One
can formulate conditions equivalent to \ref{propcor-msq} (2),
\ref{cor1-main} (3), \ref{cor2a-main-new2} (2) and
\ref{cor2a-main-new1} (1), respectively, by suppressing the
requirement that $\phi_z^{-1}(V)$ is relatively compact. This is
because $G\cdot z$ is assumed to be locally closed in $X$ and
hence every neighbourhood $V$ of $z$ contains a neighbourhood
$V_1$ of $z$ such that $\phi_z^{-1}(V_1)$ is relatively compact
(see Lemma~\ref{lem-remark}).
\end{remark}


\section{Examples}\label{sec-examples}
We apply our results in two examples. First, we calculate
upper and lower multiplicities relative to subsequences in the
transformation group described by Rieffel in \cite[Example~1.18]{rieffel}. Then we combine the ideas of \cite[Example~1.18]{rieffel} and \cite[Example~on~p.~298]{palais} to construct a
transformation group with a sequence $(G\cdot x_n)_n$ of orbits
converging to two distinct orbits $G\cdot x_0$ and $G\cdot
z_0$ such that
\begin{gather*}
M_U(\Ind\epsilon_{x_0},(\Ind\epsilon_{x_n}))
=M_L(\Ind\epsilon_{x_0},(\Ind\epsilon_{x_n}))
=2\\
M_U(\Ind\epsilon_{z_0},(\Ind\epsilon_{x_n}))=M_L(\Ind\epsilon_{z_0},(\Ind\epsilon_{x_n}))=3.
\end{gather*}
\begin{example}
\label{ex-rieffel} In \cite[Example~1.18]{rieffel}, the space $X$ is
a closed subset of $\R^3$ and the group is $G=\R$. The action is free, with all
orbits closed in $X$. The orbit space $X/G$ is a compact Hausdorff space,
and is discrete except for one limit point. This limit point is
$\lbrace (0,s,0)\colon s\in\R\rbrace$ with the action of $\R$ on it
being translation, and orbit representative $x_0=(0,0,0)$.

Let $(b_n)_{n\geq 1}$ be a strictly decreasing sequence of real
numbers converging to $0$. If $n\geq 1$ then the $n$th orbit
representative is $x_n=(b_n,0,0)$, and  assigned to each orbit is a
repetition integer $L_n\geq 0$.

If $L_n=0$ then the orbit of $x_n$ is the line $\lbrace (b_n,s,0)\colon s\in\R\rbrace$ with the action of $\R$ on it
being translation.
If $L_n\geq 1$,  choose a strictly
decreasing finite sequence $(b_n^j)_j$ of length $L_n+1$, with
$b_{n+1}< b^{L_n}_n<\cdots <b^2_n<b^1_n<b^0_n=b_n$.  Let
$x^j_n=(b^j_n,0,0)$ for $0\leq j\leq L_n$; note that $x_n=x_n^0$. The points
$x^j_n\  (j=0,\dots, L_n)$ are in the orbit $G\cdot x_n$. The $n$th
orbit consists of $L_n + 1$ vertical line segments parallel to the
$y$-axis, each of these line segments goes through  $x_n^j$
for some $j$. The vertical line segments  are joined by $L_n$ arcs.
The action is described by specifying what happens to the orbit
representatives $x_n$ $(n\geq 1)$:
\[
t\cdot x_n =
\left\{\begin{array}{ll} (b_n,t,0) &\text{if\ }t\in(-\infty,n];
\\( b_n^{L_n},s,0) &\text{if\ }s\in(-n,\infty)\text{\ and\ } t=s+L_n(2n+1);
\\( b_n^j,s,0) &\text{if\ }s\in(-n,n]\text{\ and\ } t=s+j(2n+1)\\
&\quad\quad\quad\quad\text{for\ }j\in\{1,\dots,L_n-1\};
\\\big( (1-s)b^j_n + sb^{j+1}_n, n\cos(\pi s), n\sin(\pi s) \big)
&\text{if\ }s\in(0,1]\text{\ and\ } t=s+n+j(2n+1),\\
&\quad\quad\quad\quad\text{for\ }j\in\{0,\dots,L_n-1\}.
\end{array}
\right.\]
\end{example}

\begin{lemma}\label{lem-ex}
Let $(G,X)$ be the free transformation group  from
Example~\ref{ex-rieffel} with orbit representatives $x_n$. Then
\begin{gather*}
M_L(\Ind\epsilon_{x_0},(\Ind\epsilon_{x_n}))=\liminf_n (L_n+1);
\\
M_U(\Ind\epsilon_{x_0},(\Ind\epsilon_{x_n}))=\limsup_n
(L_n+1).
\end{gather*}
\end{lemma}

\begin{proof}
Suppose that $\liminf_n (L_n+1)=k<\infty.$ We will show that $M_L(\Ind\epsilon_{x_0},(\Ind\epsilon_{x_n}))\geq k$
and $M_L(\Ind\epsilon_{x_0},(\Ind\epsilon_{x_n}))\leq k$ using Theorems~\ref{thm-main} and \ref{thm-M}.

Since $\liminf_n (L_n+1)=k$ there exists $n_0$ such that $L_n+1\geq
k$ for $n\geq n_0$. For each $n\geq n_0$ and $j\in\{0,\dots,k-1\}$ choose $t_n^{(j)}=j(2n+1)$
for . Then, as $n\to \infty$, we have
\begin{align*}
&t_n^{(0)}\cdot x_n=x_n\to x_0\\
&t_n^{(j)}\cdot x_n=(b^j_n,0,0)=q^j_n\to x_0\quad\text{for\ }j\in\{1,\dots,k-1\} \\
&t_n^{(i)}-t_n^{(j)}=(i-j)(2n+1)\to\infty\quad\text{if\ }i\neq j.
\end{align*}
So  $(x_n)_n$ converges $k$-times  in $X/G$ to $x_0$, and hence
$M_L(\Ind\epsilon_{x_0},(\Ind\epsilon_{x_n}))\geq k$ by Theorem~\ref{thm-main}.

Now consider the (relatively) open relatively compact neighbourhood
\begin{equation}\label{eq-ex1-lem}
V=\big( [0,1)\times(-1/2, 1/2)\times (-1/2, 1/2) \big)\cap X
\end{equation}
of $x_0$. Since $\liminf_n (L_n+1)=k$  there exists a subsequence
$(x_{n_i})_i$ such that $L_{n_i}+1=k$ and $x_{n_i}\in V$ for all
$i$. Note that the arcs of $G\cdot x_{n_i}$ joining the parallel
line segments of the orbit do not meet $V$.  So if $t\cdot
x_{n_i}\in V$ then $t$ is either in the interval $(-1/2,1/2)$ or in
a translate by $j(2n_i+1)$ of it, where $1\leq j\leq k-1$.  Thus
$\phi_{x_0}^{-1}(V)=(-1/2,1/2)$ is relatively compact and
\[
\nu(\phi_{x_{n_i}}^{-1}(V))
=k\nu(\phi_{x_0}^{-1}(V))<\infty
\]
for all $i$, so $M_L(\Ind\epsilon_{x_0},(\Ind\epsilon_{x_n}))\leq k$ by
Theorem~\ref{thm-M}. Thus $M_L(\Ind\epsilon_{x_0},(\Ind\epsilon_{x_n}))= k$

If $\liminf_n(L_n+1)=\infty$, then for any $k\in\P$ we have $L_n+1\geq k$ eventually,
and $(x_n)_n$ converges $k$-times to $x_0$ in $X/G$ as above.
Thus $M_L(\Ind\epsilon_{x_0},(\Ind\epsilon_{x_n}))=\infty$.

Now suppose that $\limsup_n(L_n+1)=k<\infty$. Then there exists a
subsequence $(x_{n_i})_i$ such that $L_{n_i}+1=k$ for all $i$. As
above, $(x_{n_i})$ converges $k$-times to $x_0$ in $X/G$, and hence
$M_U(\Ind\epsilon_{x_0},(\Ind\epsilon_{x_n}))\geq k$ by Theorem~\ref{thm-main-upper}.
On the other hand, there exists $n_1$ such that  for all $n\geq
n_1$, $x_n\in V$ and $L_n+1\leq k$. We have
\[
\nu(\phi_{x_n}^{-1}(V))
=k\nu(\phi_{x_0}^{-1}(V))=k<\infty
\]
whenever $n\geq n_1$, and hence $M_U(\Ind\epsilon_{x_0},(\Ind\epsilon_{x_n}))\leq k$
by Theorem~\ref{thm-Mupper}. Thus $M_U(\Ind\epsilon_{x_0},(\Ind\epsilon_{x_n}))= k$.

If $\limsup_n(L_n+1)=\infty$ then, given $k\in\P$, there exists a
subsequence $(x_{n_i})_i$ such that $L_{n_i}+1\geq k$.  Then
$(x_{n_i})$ converges $k$-times to $x_0$ in $X/G$ and
$M_U(\Ind\epsilon_{x_0},(\Ind\epsilon_{x_n}))\geq k$.  Thus
$M_U(\Ind\epsilon_{x_0},(\Ind\epsilon_{x_n}))=\infty$.
\end{proof}

\begin{remark}
In the situatation of Lemma~\ref{lem-ex} one can easily show that $M_L(\Ind\epsilon_{x_0})=M_L(\Ind\epsilon_{x_0},(\Ind\epsilon_{x_n}))$ and $M_U(\Ind\epsilon_{x_0})=M_U(\Ind\epsilon_{x_0},(\Ind\epsilon_{x_n}))$ using Lemma~\ref{lem-Lsequence} and \cite[Lemma~1.2]{AK}, respectively.
\end{remark}
\begin{example}
The idea of this example is to splice together two instances of Example~\ref{ex-rieffel},
one with constant repetition integer $1$ and the other with constant repetition integer $2$.

Again, the space $X$ is a closed subspace of $\R^3$, the action of the group $G=\R$ is free
and the orbits are closed in $X$.  The orbit space is compact but non-Hausdorff.
This time we have two limit orbits,
$\lbrace (0,s,0)\colon  s\in\R\rbrace$ and
$\lbrace (1,s,0)\colon s\in\R\rbrace$ with  orbit representatives
$x_0=(0,0,0)$  and $z_0=(1,0,0)$; the action of $\R$ on these two orbits is by translation.

If $n\geq 1$, the data for the $n$th orbit is:
\[
b_n=2^{-2n},\quad b_n^1=2^{-(2n+1)},\quad a_n=1-2^{-3n},\quad a_n^1=1-2^{-(3n+1)},\quad a_n^2=1-2^{-(3n+2)};
\]
the orbit representative is
\[x_n=\left(\frac{a_n+b_n}{2},0,n\right).\] The orbit
consists of five vertical line segments parallel to the $y$-axis,
two on the left of  and three on the right of the line
$\{(1/2,s,0):s\in\R\}$; the line segments are joined by four arcs.
Again we describe the action by specifying what happens to the orbit
representatives $x_n$. If $n\geq 1$ and  $u\in\R$ then
\[
u\cdot x_n \left\{
\begin{array}{ll} (b^1_n,u-3n-\frac{3}{2},0)
&\text{if\ }u\in[2n+\frac{3}{2},\infty);
\\\big( b_n+s(b_n^1-b_n),n\cos(\pi s),n\sin(\pi s)\big)
&\text{if\ }u\in[2n+\frac{1}{2},2n+\frac{3}{2})
\\&\quad\quad\text{and\ } s=u-2n-\frac{1}{2};
\\( b_n,u-n-\frac{1}{2},0)
&\text{if\ }u\in[\frac{1}{2}, 2n+\frac{1}{2});
\\\big(a_n+(u+\frac{1}{2})(b_n-a_n),-n\sin(u\pi),n\cos(u\pi)\big)
&\text{if\ }u\in[-\frac{1}{2},\frac{1}{2});
\\( a_n, u+n+\frac{1}{2},0 )
&\text{if\ }u\in[-2n-\frac{1}{2},-\frac{1}{2});
\\\big(a^1_n+s(a_n-a^1_n),n\cos(\pi s),n\sin(\pi s)\big)
&\text{if\ }u\in[-2n-\frac{3}{2}, -2n -\frac{1}{2})
\\&\quad\quad\text{ and\ } s=u+2n+\frac{3}{2};
\\(a^1_n, u+3n+\frac{3}{2}, 0)
&\text{if\ }u\in [-4n-\frac{3}{2},-2n-\frac{3}{2});
\\\big(a^2_n+s(a^1_n-a^2_n),n\cos(\pi s),n\sin(\pi s)\big)
&\text{if\ }u\in[-4n-\frac{5}{2}, -4n -\frac{3}{2}),
\\&\quad\quad\text{and\
 } s=u+4n+\frac{5}{2};
\\(a^2_n,u+5n+\frac{5}{2},0)
&\text{if\ }u\in(-\infty, -4n-\frac{5}{2}).
\end{array}
\right.
\]
It is straightforward to
verify that $M_U(\Ind\epsilon_{x_0},(\Ind\epsilon_{x_n}))=M_L(\Ind\epsilon_{x_0},(\Ind\epsilon_{x_n}))=2$ and
$M_U(\Ind\epsilon_{z_0},(\Ind\epsilon_{x_n}))=M_L(\Ind\epsilon_{z_0},(\Ind\epsilon_{x_n}))=3$
using the methods of Lemma~\ref{lem-ex}
\end{example}

\appendix\section{\ }\label{appendix}
The purpose of this appendix is to establish sequence versions of
\cite[Propositions 2.2 and 2.3]{AS} which were used in \S\ref{sec-measure} and \S\ref{sec-main}.  See \S\ref{sec-prelim} for the relevant notation.

\begin{lemma}\label{lemma-sepmult} Let $A$ be a separable $C^*$-algebra, $\pi\in \hat A$ and
$(\pi_n)_{n\geq 1}$ a sequence in $\hat A$.  Then there exists a
subsequence $(\pi_{n_i})_i$ such that
$$ M_L(\pi,(\pi_{n_i})_i)
=M_U(\pi,(\pi_{n_i})_i)=M_U(\pi,(\pi_n)).$$
\end{lemma}

\begin{proof} Let $\phi$ be a pure state of $A$ associated with $\pi$.
Suppose that $M_U(\pi,(\pi_n))=m$ (where possibly $m=\infty$).
It suffices to construct $(\pi_{n_i})$ such that
$M_L(\pi,(\pi_{n_i})_i)\geq m$. By assumption, there exists
$N_0\in\N$ such that, for $N\in\N$ satisfying $N\subset N_0$,
\begin{equation}\label{eq-A1}
M_U(\phi,N,(\pi_n))=\limsup_n d(\pi_n,\phi, N)=m.
\end{equation}
Since $A$ is separable, the weak$^*$-topology on $A^*$ is first
countable and so there is a basic decreasing sequence
$(N_i)_{i\geq 1}$ in $\N$ such that $N_i\subset N_0$ for all
$i\geq1$. By \eqref{eq-A1}, for all $i\geq1$,
\begin{equation}\label{eq-A2}
\limsup_n d(\pi_n,\phi, N_i)=m.
\end{equation}

Firstly, suppose that $m<\infty$.  By \eqref{eq-A2}, we may construct a
strictly increasing sequence $(n_i)_{i\geq1}$ such that
$d(\pi_{n_i},\phi,N_i)=m$ for all $i\geq1$. Let $N\in\N$. There
exists $i_0$ such that $N_{i_0}\subset N$. For $i\geq i_0$,
$$d(\pi_{n_i},\phi,N)\geq d(\pi_{n_i},\phi,N_i)=m.$$
Thus
$$M_L(\phi,N,(\pi_{n_i}))= \liminf_i d(\pi_{n_i},\phi,N)\geq m$$
and hence
$$M_L(\pi,(\pi_{n_i}))=\inf_{N\in\N}M_L(\phi,N,(\pi_{n_i}))\geq
m.$$

Secondly, suppose that $m=\infty$. By \eqref{eq-A2}, we may construct a
strictly increasing sequence $(n_i)_{i\geq1}$ such that
$d(\pi_{n_i},\phi,N_i)\geq i$ for all $i\geq1$. Let $N\in\N$.
There exists $i_0$ such that $N_{i_0}\subset N$. For $i\geq
i_0$,
$$d(\pi_{n_i},\phi,N)\geq d(\pi_{n_i},\phi,N_i)\geq i.$$
Thus
$$M_L(\phi,N,(\pi_{n_i}))= \liminf_i d(\pi_{n_i},\phi,N)=\infty$$
and hence
$$M_L(\pi,(\pi_{n_i}))=\inf_{N\in\N}M_L(\phi,N,(\pi_{n_i}))=\infty.$$
\end{proof}


\begin{lemma}\label{lem-Lsequence}
Let $A$ be a separable $C^*$-algebra.  Let $\pi\in\hat A$ such that $\{\pi\}$ is not open in $\hat A$.  Then there exists a sequence $(\pi_n)_{n\geq 1}$ in $\hat A$ such that $\pi_n\neq \pi$ for all $n\geq 1$, $\pi_n\to_n\pi$,  and
\[M_L(\pi)=M_L(\pi,(\pi_n))=M_U(\pi,(\pi_n)).\]
\end{lemma}

\begin{proof}
First assume that $M_L(\pi)=\infty$.  Since $\hat A$ is second countable and $\{\pi\}$ is not open in $\hat A$, there exists a sequence $(\pi_n)_{n\geq 1}$ in $\hat A$ such that $\pi_n\to \pi$ and $\pi_n\neq \pi$ for all $n$. Then
\[\infty=M_L(\pi)\leq M_L(\pi,(\pi_n))\leq M_U(\pi,(\pi_n)),\]
and we must have equality throughout.

Next, assume that $M_L(\pi)=k$ where $1\leq k<\infty$.  Let
$\phi$ be a pure state associated with $\pi$.  There exists
$N_0\in\N$ such that if $N\in\N$ and $N\subset N_0$, then
$M_L(\phi, N)=\liminf_{\sigma\to\pi,\sigma\neq\pi}
d(\sigma,\phi,N)=k.$ Since $A$ is separable, there exists a
decreasing base $(V_n)_{n\geq 1}$ of open neighbourhoods of $\pi$
in $\hat A$.  We have
\[
\liminf_{\sigma\to\pi,\sigma\neq\pi} d(\sigma,\phi,N_0)=k,
\]
so for $n\geq 1$ there exist $\pi_n\in V_n\setminus\{\pi\}$ such that $d(\pi_n,\phi, N_0)=k$.
Since $\{V_n\}$ is decreasing, $\pi_n\to \pi$ as $n\to\infty$.  Note that
\begin{align*}
M_U(\pi,(\pi_n))&=\inf_{N\in\N}M_U(\phi, N,(\pi_n))\\
&\leq M_U(\phi, N_0,(\pi_n))\\
&=\limsup_n d(\pi_n,\phi, N_0)=k.
\end{align*}
Now $k=M_L(\pi)\leq M_L(\pi,\pi_n))\leq M_U(\pi,(\pi_n))\leq k$ and the result follows.
\end{proof}


\end{document}